\documentclass[11pt]{article}

\usepackage[centertags]{amsmath}
\usepackage{amsfonts}
\usepackage{amssymb}
\usepackage{theorem}
\usepackage{ulem}
\usepackage{isolatin1}
\usepackage{epsfig}
\usepackage{subfigure}
\usepackage[all]{xy}
\usepackage{changebar}

\theoremstyle{break} \newtheorem{theorem}{Theorem}[section]
\theoremstyle{break} 
\theoremstyle{break} \newtheorem{definition}{Definition}[section]       
\theoremstyle{break} \newtheorem{lemma}[theorem]{Lemma}
\theoremstyle{break} 
\theoremstyle{break} \newtheorem{example}[theorem]{Example}
\theoremstyle{break} 
\theoremstyle{break}
{\theorembodyfont{\rmfamily}\newtheorem{remark}[theorem]{Remark}}
{\theorembodyfont{\rmfamily}}
\theoremstyle{break} 
\theoremstyle{break} 
\theoremstyle{break} \newtheorem{proposition}[theorem]{Proposition}
\theoremstyle{break} 
\numberwithin{equation}{section}

\newcommand{\R}{{\mathbb{R}}}
\newcommand{\D}{{\mathbb{D}}}
\newcommand{\C}{{\mathbb{C}}}

\def\d{\mathop{{\rm d}_{\lambda}}}

\def\Re{\mathop{{\rm Re}}}

\begin{document}

\renewcommand{\thefootnote}{}
\stepcounter{footnote}
\begin{center}
{\bf \Large 
The behaviour of solutions of the Gaussian curvature\\[2mm] equation near an isolated boundary point
\footnote{2000 Mathematics Subject
 Classification: Primary  35J60, 32F45, 53A30\\
The first author was supported by a HWP scholarship; the second author
received partial support  from the German--Israeli Foundation (grant G--809--234.6/2003).}}
\end{center}

\medskip

\renewcommand{\thefootnote}{\arabic{footnote}}
\setcounter{footnote}{0}
\begin{center}
{\large  \bf Daniela Kraus and Oliver Roth}\\[2mm]
{\small Universit\"at W\"urzburg,
 Mathematisches Institut, \\[-0.5mm]     
 D--97074 W\"urzburg, Germany\\[1mm]
dakraus@mathematik.uni-wuerzburg.de\\ roth@mathematik.uni-wuerzburg.de}\\[1mm]
\end{center}

\centerline{\today}

\renewcommand{\thefootnote}{\arabic{footnote}}

\medskip
\begin{center}
\begin{minipage}{13cm} {\bf Abstract.}   
A classical result of Nitsche \cite{Nit57} about the
behaviour of the solutions to the Liouville equation $\Delta u=4
\, e^{2u}$ near isolated  singularities is generalized  
to solutions of the Gaussian curvature equation
 $\Delta u=- \kappa(z)\,  e^{2u}$ where $\kappa$ is
a negative H\"older continuous function. As an application
 a higher--order version of the Yau--Ahlfors--Schwarz lemma
for complete conformal Riemannian metrics is obtained.
\end{minipage}
\end{center}

\section{Introduction}
In \cite{Nit57} Nitsche gave a detailed description of the behaviour of 
 the real--valued solutions to the Liouville equation
\begin{equation}\label{eq:pde1}
\Delta u=4 \, e^{2u} \,  
\end{equation}
on plane domains near their isolated singularities. The first purpose of the present
paper is to extend Nitsche's results to the solutions of the more general Gaussian
curvature equation
\begin{equation}\label{eq:pde2}
\Delta u= -\kappa (z)\,  e^{2u}
\end{equation}
with strictly negative H\"older continuous
functions $\kappa(z)$. 
It suffices to consider this PDE on the
open unit disk $\D=\{ z \in \C\, : \, |z|<1\}$ in the complex plane $\C$.
We use the notation
$$ M_u(r):=\sup \limits_{|z|=r} u(z) \, $$
for real--valued functions $u$ defined in a punctured neighborhood of $z=0$
and call 
$$ \alpha(u):=\lim \limits_{r \searrow 0} \frac{M_u(r)}{\log(1/r)}$$
the order of $u$ if this limit exists.

\begin{theorem}\label{thm:hoelder}
Let $\kappa:\D \to \R$  be a  locally H\"older
continuous function with $\kappa(0)<0$. If $u:
\D \backslash \{ 0 \} \to \R$ is a $C^2$--solution to
$\Delta u= -\kappa (z)\,  e^{2u}$ in  $\D\backslash \{ 0 \}$, then 
$u$ has order $\alpha \in (-\infty,1]$ and
\begin{alignat}{3}
u(z)&= -\alpha \, \log{|z|} + v(z)\, , \quad \quad \qquad & \text{if }\alpha <1\, , \label{eq:hoelder}\\[2mm]
u(z)&= - \log|z|-\log{\log{(1/|z|)}}+w(z)\, ,\quad \quad \qquad & \text{if } \alpha=1\, ,  \label{eq:hoelder2}
\end{alignat}
where the remainder functions $v$ and $w$ are continuous in $\D$. Moreover, the first
partial derivatives with respect to $z$ and $\overline{z}$,
\begin{alignat*}{3}
& v_z(z), v_{\overline{z}}(z) \quad \text{are continuous at}\, \, z=0\,\,
 \, \quad &\text{if}  \quad &\alpha <1/2\, ;
\\[-4mm]
\intertext{and}\\[-10mm]
&  v_z(z), v_{\overline z}(z)=O(1) &\text{if} \quad &\alpha=1/2\, ,\\[1mm]
&v_z(z), v_{\overline z}(z)=O\left(|z|^{1-2\alpha}\right)  
&\text{if} \quad  &1/2<\alpha<1\, ,\\[1mm]
&w_z(z), w_{\overline{z}}(z)=O\left(|z|^{-1}\, (\log{(1/|z|)})^{-2}\right) \qquad &\text{if}
 \quad  & \alpha=1 \, ,
\end{alignat*}
when $z$ approaches $z=0$. In addition, the second partial derivatives,
\begin{alignat*}{3}
&\displaystyle v_{z\,\! z}(z), v_{z \overline{z}}(z) \text{ and } 
v_{\overline{z}\, \overline{z}}(z) \quad \text{are continuous
 at}\, \, z=0\, \,
 \, \quad &\text{if}  \quad &\alpha \le 0\, ;
\\[-4mm]
\intertext{and}\\[-10mm]
&v_{z\,\!z}(z), v_{z \overline{z}}(z),
v_{\overline{z}\,\overline{z}}(z)=O\left(|z|^{- 2\alpha}\right)
  &\text{if} \quad  &0<\alpha<1\, ,\\[1mm]
&w_{z\, \!z}(z), w_{z \overline{z}}(z), w_{\overline{z}\,\overline{z}}(z)=O\left(|z|^{-2}\,
 (\log{(1/|z|)})^{-2}\right) \qquad   &\text{if} \quad  & \alpha=1\, ,
\end{alignat*}
when $z$ tends to $z=0$.
\end{theorem}

Theorem \ref{thm:hoelder} merits some comment.
Firstly, the special case $\kappa(z) \equiv -4$ of Theorem \ref{thm:hoelder}
is Nitsche's theorem; see \cite[Satz 1]{Nit57}.\footnote{Nitsche
  considers the PDE $\Delta U=e^{U}$, which is obtained from $\Delta u=4 \, e^{2
    u}$ using the transformation $U(z)=2 \,  u(z)+\log 8$; {\it cf}.~also Remark
  \ref{rem:nitsche} below.}
Nitsche's proof is based on an ingenious application of
Liouville's classical representation formula \cite{Lio1853}
for the solutions to (\ref{eq:pde1}). Roughly speaking,
Liouville's result  says that in any disk $D \subseteq
\D \backslash \{ 0\}$ every solution $u$ to (\ref{eq:pde1}) can be written as
$$u(z)=\log \frac{|f'(z)|}{1-|f(z)|^2}\, ,$$
where $f$ is a holomorphic function in $D$. 
Now a careful and clever study of the analytic continuation of  $f$ along a path 
surrounding the singularity $z=0$ enables Nitsche to prove Theorem
\ref{thm:hoelder} for $\kappa(z) \equiv -4$. The same argument was later 
used by Yamada \cite{Yam88} and also by
Chou and Wan \cite{CW94,CW95} who were apparently unaware of Nitsche's paper.
Clearly, Nitsche's method cannot be applied to prove Theorem \ref{thm:hoelder}
for non--constant functions $\kappa(z)$ as there is no representation formula
of Liouville--type in this case.

\medskip

Secondly, the motivation for extending Nitsche's result  mainly  comes
 from the geometric interpretation of the PDE (\ref{eq:pde2}):
every solution $u$ to (\ref{eq:pde2}) gives rise to a conformal Riemannian metric
$e^{u(z)} \, |dz|$ with Gaussian curvature $\kappa(z)$ and vice versa; see Paragraph
\ref{par:conformal} below.
Thus passing from Liouville's equation (\ref{eq:pde1})
to the more general equation (\ref{eq:pde2}) 
 amounts to passing from constantly curved conformal Riemannian
metrics to metrics with variable curvature.

\medskip

We note that the constant curvature case is intimately
connected with the uniformization problem for Riemann surfaces and
the classical Schwarz--Picard problem. Its study was pursued by
Schwarz \cite{Sch1891}, Picard \cite{Pic1893,Pic1905}, Poincar\'e
\cite{Poi1898}, Bieberbach \cite{Bie12,Bie16}, Heins \cite{Hei62} and many others.
It has only recently  led to a complete proof of the uniformization theorem
for Riemann surfaces by Mazzeo and Taylor \cite{MT2002} solely based on curvature
considerations. The Schwarz--Picard problem is the problem
of investigating the solutions to $\Delta u=4 \, e^{2 u}$ {\it with prescribed
singularities} on a compact Riemann surface. The solution
to the existence--and--uniqueness part of the Schwarz--Picard
problem is due to Heins \cite{Hei62}, while the sharp growth and
regularity properties of the corresponding solutions at the singularities are described by
Nitsche's theorem.

\medskip

The more general and difficult case of variable curvature is strongly related to the  Berger--Nirenberg problem in
differential geometry (see Aviles \& McOwen \cite{AM85}, Troyanov 
\cite{Tro90}, Hulin \& Troyanov \cite{HT92}, and Chang \cite{Chang} as some
of the many references). In particular, in analogy to the Schwarz--Picard
problem the solutions to the equation
$\Delta u=-\kappa(z) \, e^{2 u}$ with prescribed singularities
were studied in \cite{Tro90,HT92,McO93} and  
 existence and uniqueness of solutions
for strictly negative curvature functions were obtained. Theorem \ref{thm:hoelder}
supplements these results  by
extending Nitsche's theorem to the   variable curvature case and
thus establishing the corresponding sharp growth and regularity properties of the
solutions near their singularities.

\medskip


The basic ingredients we employ in the present paper to carry
over all  of Nitsche's  results for the constant curvature case to  
the variable curvature case
are a generalized maximum principle for the Gaussian curvature equation
({\it cf}.~Theorem \ref{thm:ext_max} below), which allows an application of 
the method  of sub-- and supersolutions, and potential--theoretic tools. Our approach 
reveals precisely how the growth and regularity of the
remainder
functions $v$ and $w$
at $z=0$ depend on the regularity of the curvature function $\kappa$. It
also leads to a number of refinements of Theorem \ref{thm:hoelder}
with weaker assumptions. These refinements will be discussed in detail in
Section \ref{sec:results}.
For instance, we shall see that if 
$u(z)$ is of the form (\ref{eq:hoelder2}), then
the function $w$ is continuous,  when $\kappa(z)$ is only assumed to be continuous 
at $z=0$ ({\it cf}.~Theorem \ref{thm:alpha=1,continuity} and Example \ref{ex:alpha=1,kappa=beschraenkt}).

\medskip

As indicated above Theorem \ref{thm:hoelder} and its refinements (see Section \ref{sec:results})
give precise information about the behaviour of regular\footnote{We call a
  conformal metric 
$\lambda(z)\,|dz|$ on a domain $\Omega \subseteq \C$ regular, if 
its density $\lambda$ is of class $C^2$ in $\Omega$.}
 conformal Riemannian metrics
and their first and second derivatives near isolated singularities.
In order to state these information in more geometric terms
we first recall the definition of two natural
derivatives associated with regular conformal metrics. The connection (or
Pre--Schwarzian or Christoffel symbol) of a regular conformal metric $\lambda(z) \,
|dz|$ on a plane domain $\Omega \subset \C$ is defined by
$$ \Gamma_{\lambda}(z)=2 \frac{\partial \log \lambda(z)}{\partial z} $$
and the Schwarzian of $\lambda(z) \, |dz|$ is given by
$$ S_{\lambda}(z)=\frac{\partial \Gamma_{\lambda}(z)}{\partial z}-\frac{1}{2}
\Gamma_{\lambda}(z)^2= 2 \left[  \frac{\partial^2 \log \lambda(z)}{\partial
    z^2}- \left(  \frac{\partial \log \lambda(z)}{\partial z} \right)^2
\right] \, .$$
These differential quantities obey simple transformation laws under conformal
change of coordinates; see \cite{Min97} and \cite{Sch}.

\begin{theorem} \label{thm:2}
Let $\lambda(z) \, |dz|$ be a regular conformal Riemannian metric on a  domain $\Omega \subset
\C$ with an isolated boundary point at $z=0$, and suppose that the curvature
 $\kappa : \Omega \to \mathbb{R}$ has a H\"older continuous extension to 
$\Omega \cup \{ 0 \}$ such that $\kappa(0)<0$. Then $\log \lambda$ has order
$\alpha \in (-\infty,1]$ and
\begin{itemize}
\item[(a)]
 $\displaystyle \lim \limits_{z \to 0} \left( |z| \log(1/|z|) \right)
 \lambda(z)=\begin{cases} 0 & \text{ if } \alpha <1 \\
   1/\sqrt{-\kappa(0)} & \text{ if } \alpha=1\, ;\end{cases}$
\item[(b)]
 $\displaystyle 
\lim \limits_{z \to 0}  z \, \Gamma_{\lambda}(z)=-\alpha$;
\item[(c)] $\displaystyle 
\lim \limits_{z \to 0}  z^2 \, S_{\lambda}(z)=\alpha \, (2-\alpha)/2$.
\end{itemize}
\end{theorem}

\begin{remark}
Theorem \ref{thm:2} extends a result of D.~Minda \cite{Min97} who proved
Theorem \ref{thm:2} for the special case that $\lambda(z) \, |dz|$ is the
hyperbolic metric in $\Omega$ (with constant curvature
$-4$). He asked (\cite[\S 4]{Min97}) for an generalization
of his results to metrics with {\it constant} negative curvature. Theorem
\ref{thm:2} gives a complete answer to Minda's question even for the much more
general case of metrics with {\it variable} negative curvature.
\end{remark}

Another application of Theorem \ref{thm:hoelder} deals with a higher--order
version of the so--called Yau--Ahlfors--Schwarz lemma; see Yau \cite{Yau78}
and Ahlfors \cite{Ahl38}.
In the complex one--dimensional case 
Yau's generalized Ahlfors--Schwarz lemma  says that if 
$\lambda_{\Omega} \, |dz|$ is the hyperbolic metric   of a (hyperbolic) domain $\Omega \subset \C$ and $\lambda(z) \,
|dz|$ is a {\it complete} regular conformal metric in $\Omega$ with curvature
$\kappa(z) \ge -4$, then $\lambda(z) \ge \lambda_{\Omega}(z)$ for all $z \in \Omega$.
Thus Yau's lemma derives a global estimate for conformal metrics 
from global assumptions. There are boundary versions
of these results ({\it cf}.~Bland \cite{Bla83}, Troyanov \cite{Tro1992}, and Kraus,
Roth \& Ruscheweyh \cite{KRR06}). The following theorem gives  precise 
{\it local} information about  conformal metrics {\it and} its derivatives up to
second order at isolated boundary points under {\it local} assumptions.
We call a conformal metric $\lambda(z) \, |dz|$ on a domain $\Omega$ with an
isolated boundary point at $z=0$ 
{\it locally complete} (at $z=0$) if
$$ \lim \limits_{\overset{z\to 0}{z \in \Omega}} \d(z_0,z) = +\infty$$ 
for some (and then every) point $z_0 \in \Omega$. Here $\d$ denotes the distance function induced
by $\lambda(z) \, |dz|$. A basic reference about locally complete metrics 
is the work of Huber \cite{Hub}; see also \cite{Bla83} and \cite{KRR06}.

\begin{theorem}[Higher--order Yau--Ahlfors--Schwarz Lemma] \label{thm:3}
Let $\Omega$ be a hyperbolic domain  with an 
isolated boundary point at $z=0$ and let $\lambda_{\Omega}(z) \, |dz|$ be the
hyperbolic metric of $\Omega$ (with curvature $-4$).
Let $\lambda(z) \, |dz|$ be a locally complete conformal Riemannian metric on $\Omega$
 with H\"older continuous curvature $\kappa : \Omega \cup \{0\} \to
 \mathbb{R}$ such that $\kappa_{\lambda}(0)=-4$. Then
\begin{itemize}
\item[(a)]
 $\displaystyle \lim \limits_{z \to 0} \frac{\lambda(z)}{\lambda_{\Omega}(z)}=1$
\item[(b)]
 $\displaystyle 
\lim \limits_{z \to 0}  \frac{\Gamma_{\lambda}(z)}{\Gamma_{\lambda_{\Omega}}(z)}=1 $
\item[(c)] $\displaystyle 
\lim \limits_{z \to 0}  \frac{S_{\lambda}(z)}{S_{\lambda_{\Omega}}(z)}=1$.
\end{itemize}
\end{theorem}

The paper is organized as follows. In
 Section \ref{sec:preliminaries} we begin with a discussion of
some basic facts from conformal geometry (\S \ref{par:conformal}) and
prove the extended maximum principle for the curvature equation
(\S \ref{par:max}).  
Section \ref{sec:results} contains the statements and  proofs of 
the main results of this work. It starts in Paragraph \S \ref{par:conti}
with a proof of the representation
formulas (\ref{eq:hoelder}) and (\ref{eq:hoelder2}) of Theorem \ref{thm:hoelder}
 under minimal hypotheses
on the curvature function (Theorem \ref{thm:beschraenkt})
and a discussion, when the remainder functions $v(z)$ and $w(z)$ are
continuous
(see Theorem \ref{thm:alpha=1,continuity}).
We then establish the growth and regularity properties
of the first  derivatives of $v(z)$
and $w(z)$ in \S \ref{par:abl1}--\S \ref{par:proofs} 
and of the second derivatives in  \S \ref{par:abl2} 
again under minimal
assumptions on the curvature function. Theorems \ref{thm:2} and \ref{thm:3}
are proved in \S \ref{par:thm23}.
All of our results are essentially
sharp as  will be illustrated with a number of examples. The paper ends with
an appendix on  the required (non--standard) tools from potential theory.

\section{Preliminaries}\label{sec:preliminaries}

As indicated above the proof of Theorem \ref{thm:hoelder} uses a mixture  of
different methods from conformal geometry, subharmonic functions including a
maximum--principle, and potential theory.  In this preparatory section, we 
first recall some basic facts about conformal Riemannian metrics and show that Theorem
\ref{thm:hoelder} is best possible (\S \ref{par:conformal}). 
Paragraph \ref{par:max} is devoted to an
 extended maximum principle for the Gaussian curvature equation which is
  the main technical tool for the proof of
Theorem \ref{thm:hoelder}.

\subsection{Conformal Riemannian metrics and the Gaussian curvature equation}
\label{par:conformal}

Every positive upper semi--continuous function
$\lambda$ on a domain $G \subset \C$ induces a conformal Riemannian metric
 $\lambda(z)\, |dz|$ on $G$. 
The (Gaussian) curvature of a regular conformal Riemannian metric $\lambda(z)\,|dz|$
is defined  by
\begin{equation*}
\kappa_{\lambda}(z):=-\frac{(\Delta \log{\lambda})\, (z)}{\lambda(z)^2} \, ,
\qquad \text{ where } \Delta=\frac{\partial^2}{\partial
  x^2}+\frac{\partial^2}{\partial y^2} \, , \quad z=x+i y\, .  
\end{equation*}
Thus, if $\lambda(z)\, |dz|$ is a regular conformal Riemannian metric on a domain $G$
    with curvature $\kappa_{\lambda}$, then the
    function $u:=\log \lambda$ is a $C^2$--solution to the  Gaussian curvature
    equation  
$$\Delta u= -\kappa_{\lambda}(z) \, e^{2u}$$
in $G$ and vice versa.

\medskip

A basic property of curvature is its absolute conformal invariance. This
means that the pullback of a conformal Riemannian metric $\lambda(z)\, |dz|$ on a
domain $D$, defined by $(f^*\lambda)(z) \, |dz|:=
\lambda(w)\, |dw|$, where $w=f(z)$ is a holomorphic map from  a domain $G$ to
$D$, is  a conformal Riemannian metric on $G$ off the set of critical points of $f$ 
with curvature
$$\kappa_{f^*\lambda}(z)=\kappa_\lambda(f(z))\, .$$  
This conformal invariance provides a simple, but flexible tool to construct
new conformal Riemannian metrics from old ones. For instance, the Poincar\'e metric
$$ \frac{|dz|}{2 \, |z| \, \log(1/|z|)}$$
on the punctured  unit disk $\D\backslash \{ 0 \}$ 
has constant curvature $-4$. Pulling back this metric via the map
$f(z)=z/R$ for some $R>1$ gives another conformal Riemannian metric on $\D \backslash \{
0 \}$ with curvature $-4$.
The following examples have been constructed  along these lines.

\begin{example} \label{ex:dani1}
The function
\begin{equation*}
u_{\alpha}(z)= 
\begin{cases}
\displaystyle - \alpha \log |z|+ \log \left( (1-\alpha)\,  \frac{4^{2-\alpha}}{2}\,    \frac{ \, |1+z|\,
 \displaystyle{ \left|2+ z
   \right|^{-\alpha}}}{4^{2(1-\alpha)}-|2z+z^2|^{2(1-\alpha)}}\right)\quad
 &\text{if} \, \, \, \, \alpha \le 0 \, ,\\[4mm]
\displaystyle - \alpha \log |z| + \log \left( \frac{1-\alpha}{1-|z|^{2(1-\alpha)}} \right) \quad
 &\text{if} \, \, \, \, 0 <\alpha < 1 \, ,  \\[4mm]
\displaystyle- \log|z| - \log \log \frac {1}{|z|}  +
 \log\left(\frac{1}{2} \cdot \frac{\log(1/|z|)}{1 +
\log(1/|z|)}\right) &\text{if}  \, \,\, \,  \alpha=1\, .
\end{cases}
\end{equation*}
is a $C^2$--solution to $\Delta u =4\, e^{2u}$ in $\D \backslash \{ 0 \}$.
Thus all cases $\alpha \le 1$ in Theorem
\ref{thm:hoelder} do occur. In the notation of Theorem \ref{thm:hoelder}
we obtain for the partial derivatives of the remainder functions $v$ (if $\alpha<1$)
and $w$ (if $\alpha=1$)  
\begin{alignat*}{2}
\displaystyle &\lim_{z \to 0} v_z(z)=\frac{1}{2} - \frac{\alpha}{4} \qquad 
&\text{if}\, \, \, 
&\alpha \le 0\, ,\\[2mm]
\displaystyle &v_z(z)=\frac{1-\alpha}{1-|z|^{2(1-\alpha)}}
\frac{\overline{z}}{|z|^{2\alpha}}&\text{if}\, \, \, &
0< \alpha < 1\, ,\\[2mm]
&\displaystyle w_z(z)=-\frac{1}{2 \, z \, \left(\log(1/|z|)\right)^2} \left(
  \frac{ \log(1/|z|)}{1 + \log(1/|z|)} \right)&\text{if}\, \, \, &\alpha=1 \, ,\\[-3mm]
\intertext{and} \\[-10mm]
&\displaystyle \lim_{z \to 0}v_{zz}(z)=-\frac{1}{2} + \frac{\alpha}{8} \qquad
&\text{if}\, \, \,
&\alpha \le 0\, ,\\[2mm]
&\displaystyle  v_{zz}(z)=(1-\alpha) \frac{\overline{z}}{z}
\frac{|z|^{2(1-\alpha)}-\alpha}{\left(1-|z|^{2(1-\alpha)}\right)^2}
  \frac{1}{|z|^{2 \alpha}}
\quad &\text{if}\, \, \,  &0<\alpha <1\, , \\[4mm]
&\displaystyle w_{zz}(z)= - \frac{1}{4\, z^2 \left( \log(1/|z|)
  \right)^2}  
\frac{ \displaystyle\frac{1}{(\log(1/|z|))^2}-2}{\left(1+\displaystyle 
\frac{1}{\log(1/|z|)} \right)^2}
\quad &  \quad \text{if}\, \, \,  & \alpha=1\, .
\end{alignat*}
\end{example}

In particular, all the statements of Theorem \ref{thm:hoelder} are best
possible even  for the constant curvature case.

\begin{remark} \label{rem:nitsche}
We note at this point that Nitsche \cite[Satz 1]{Nit57} erroneously asserts 
that the first partial derivatives $v_z, v_{\overline{z}} = O(|z|^{1-2
  \alpha})$ and the second partial
derivatives of $v$ are in $O(|z|^{-2 \alpha})$ for all $\alpha <1$. This, however, as
the above examples show, is only true for $1/2 < \alpha<1$ and $0< \alpha<1$
respectively.
\end{remark}

For completeness we also notice that in Theorem \ref{thm:hoelder} the condition
 that $\kappa(z)$ is bounded from above
and below by negative constants at least close to $z=0$ cannot be dropped
completely. For instance,
the function
$$u(z)=\frac{1}{2}\log \left(  \frac{1}{|z|\, \log(e/|z|)} \right)=
-\frac{1}{2} \log |z|-\frac{1}{2} \log \log \frac{1}{|z|}
+\frac{1}{2} \log \left( \frac{\log (1/|z|)}{1+\log( 1/|z|)} \right) $$
is a $C^2$--solution  to the PDE
$$\Delta u=\frac{1}{|z|}\, \frac{1}{2-2\log|z|}\, e^{2u}\, $$
in $\D \backslash \{ 0 \}$.
In this case $\kappa(z) \to -\infty$ as $z \to 0$.

\subsection{Subharmonic functions and an extended maximum principle} \label{par:max}

If the curvature function $\kappa(z)$ is non--positive, then every solution
to $\Delta u=-\kappa(z) \, e^{2 u}$
is obviously subharmonic. In order to exploit this  property we need to make
use of a number of facts about subharmonic functions.

\medskip

For convenience we shall reserve the notation 
$K_R:=K_R(0)$ for the open disk with  center $0$ and
radius $R$.
Let $u$ be a subharmonic function on the punctured disk $K_R\backslash
\{ 0 \}$ with $u \not \equiv -\infty$. For $0 <r < R$ let
$$ M_u(r):=\sup \limits_{|z|=r} u(z) \, .$$
Then  $M_u(r)$ is a convex function of $\log r$ (see \cite[p.~67--68]{Hay76}), so
the left and right derivatives of $r \mapsto M_u(r)$ exist everywhere in
$0<r<R$ and are equal outside a countable set. We denote the derivative by
$M_u'(r)$. Also, $r M_u'(r)$ is
monotonically increasing and
\begin{equation} \label{eq:sub1}
 \lim \limits_{r \to 0} r M_u'(r)=- \lim \limits_{r \to 0}
\frac{M_u(r)}{\log (1/r)} \in [-\infty, + \infty)
\end{equation}
exists ({\it cf}.~\cite[p.~67/68]{Hay76}). In particular, if 
\begin{equation} \label{eq:sub2}
\lim \limits_{r \to 0} \frac{M_u(r)}{\log (1/r)}=0 \, , 
\end{equation}
then $M_u'(r) \ge 0$ for $0<r <R$, so $M_u(r)$ is monotonically increasing,
and $u$ is therefore bounded above on $\overline{K_r} \backslash \{ 0 \}$ for $r<R$.
But then $u$ has a subharmonic extension to $K_R$, when we put
$$ u(0):=\limsup_{z \to 0} u(z) \in [-\infty, +\infty) \, , $$
see \cite[p.~48]{Rado}. Thus (\ref{eq:sub2}) implies that $z=0$ is a removable
singularity of the subharmonic function $u$. We will need the following
variant of this simple fact:

\begin{lemma} \label{lem:sub1}
Let $u$ be an upper semi--continuous function and let $v$ be a subharmonic function on $K_R \backslash \{ 0
\}$ with $u,v \not \equiv -\infty$ such that
\begin{itemize}
\item[(i)] $\displaystyle \lim \limits_{r \to 0}
  \frac{M_u(r)}{\log(1/r)}=\lim \limits_{r \to 0}
  \frac{M_v(r)}{\log(1/r)}<+\infty$; and
\item[(ii)] $u-v$ is a non--negative subharmonic function in $K_R
  \backslash \{ 0 \}$.
\end{itemize}
Then $u-v$ is subharmonic in $z=0$.
\end{lemma}

{\bf Proof.}\\ Clearly, $u=(u-v)+v$ is subharmonic in $K_R \backslash
\{ 0 \}$. 
Denoting the common value of the limits in condition (i) by $\alpha \in \R$,
we observe that
$$ \lim \limits_{r \to 0} \frac{M_{u_{\alpha}}(r)}{\log(1/r)}=
\lim \limits_{r \to 0} \frac{M_{v_{\alpha}}(r)}{\log(1/r)}=0$$ 
for $u_{\alpha}(z)=
u(z)-\alpha \log (1/|z|)$ and $v_{\alpha}(z)=v(z)-\alpha \log (1/|z|)$.
Clearly, $u_{\alpha}$ and $v_{\alpha}$ are subharmonic in $K_R \backslash \{ 0 \}$.
By the discussion preceding Lemma \ref{lem:sub1}, we thus obtain 
that $u_{\alpha}$ and $v_{\alpha}$ are
subharmonic on the entire disk $K_R$, so by  \cite[p.~48 and p.~78]{Ransford}
\begin{equation*}
\begin{split}
 \lim \limits_{r \to 0} \frac{1}{2 \pi} \int \limits_{0}^{2 \pi}\frac{ 
 u(re^{it}) -v(re^{it})}{\log(1/r)}\, dt
&=\lim \limits_{r \to 0 }  \frac{1}{2 \pi} \int \limits_{0}^{2 \pi}\frac{ 
 u_{\alpha}(re^{it})}{\log(1/r)} \, dt - \lim \limits_{r \to 0 } \frac{1}{2 \pi} \int \limits_{0}^{2
 \pi}\frac{v_{\alpha}(re^{it})}{\log (1/r)}\, dt\\[2mm]
&=\lim_{r\to 0} \frac{M_{u_{\alpha}}(r)}{\log(1/r)}-\lim_{r\to 0}
 \frac{M_{v_{\alpha}}(r)}{\log(1/r)}=0 \, . 
\end{split}
\end{equation*}
Hence a result of Brelot and Saks (\cite[p.~49]{Rado}) guarantees that the
non--negative function $u-v$ is subharmonic at $z=0$, as required.\hfill{$\blacksquare$}

\bigskip

Lemma \ref{lem:sub1} is now used to 
derive a maximum principle for subharmonic sub-- and supersolutions to the PDE
$\Delta u= -\kappa(z) \, e^{2u}$ on the punctured disk $K_R \backslash\{ 0 \}$.
  The only information we need at $z=0$ is encoded
in the value of the limit (\ref{eq:sub1}). 
This  is indeed  a very useful gain in flexibility and as
we shall see later the key to Theorem \ref{thm:hoelder}.
Recall that a $C^2$--function $u$ is a subsolution to $\Delta u= -\kappa(z) \,
e^{2u}$ 
if $\Delta u \ge -\kappa(z) \, 
    e^{2u}$ and a supersolution to $\Delta u= -\kappa(z)\,  e^{2u}$ 
if $\Delta u \le -\kappa(z) \,  e^{2u}$.

\begin{theorem}[Extended maximum principle for the curvature equation]\label{thm:ext_max}
Let $\kappa : K_R \backslash \{ 0 \} \to \R$ be a non--positive function
 and let $u_1, u_2 : K _R\backslash \{ 0 \}
\to \R$, where
\begin{itemize}
\item[(i)] $u_2$ is a subharmonic supersolution to $\Delta u=-\kappa(z) \, e^{2
    u}$ in $K_R \backslash \{ 0 \}$;
\item[(ii)] $u_1$ is a subsolution to $\Delta u=-\kappa(z) \, e^{2
    u}$ in $K_R \backslash \{ 0 \}$;
\item[(iii)]  $\limsup_{z\to \xi} u_1(z) \le \liminf_{z \to \xi} u_2(z)$ for
  every $\xi\in \partial K_R$; and
  \item[(iv)] $\displaystyle \lim \limits_{r \to 0}
  \frac{M_{u_1}(r)}{\log(1/r)}\le\lim \limits_{r \to 0}
  \frac{M_{u_2}(r)}{\log(1/r)}<+\infty$. 
\end{itemize}
Then $u_1 \le u_2$ in $K_R \backslash  \{ 0 \}$.
\end{theorem}

\begin{remark}
Theorem \ref{thm:ext_max} can be easily extended to much more general
nonlinear (twodimensional) PDEs such as  the class of PDEs discussed in \cite[\S
2.3]{Jo} and \cite[\S 10.1]{GT97}. 
\end{remark}

{\bf Proof.}\\
In order to prove Theorem \ref{thm:ext_max} we apply Lemma \ref{lem:sub1}. At first  we
define on $K_R\backslash \{ 0 \}$ the subharmonic function 
$$w_1(z)=\max\{ u_1(z), u_2(z)\}\, .$$
Then the function
$$w_2(z):=w_1(z) -u_2(z)$$ 
is non--negative and  by  condition (i) and  (ii) subharmonic on
$K_R\backslash\{0\}$. In fact, if $w_2(z_0)>0$ at some point $z_0 \in K_R
\backslash \{ 0 \}$, then $w_2(z)=u_1(z)-u_2(z)>0$ in a neighborhood of
$z_0$. Thus
$$ \Delta w_2(z)=\Delta u_1(z)-\Delta u_2(z) \ge -\kappa(z) \left( e^{2
    u_1(z)}-e^{2 u_2(z)} \right) \ge 0$$
there, i.e., $w_2$ is subharmonic in this neighborhood. If $w_2(z_0)=0$ for
some point $z_0 \in K_R \backslash \{ 0 \}$, then $w_2$ satisfies the submean
inequality
$$ w_2(z_0) =0 \le \frac{1}{2\pi} \int \limits_{0}^{2 \pi} w_2\left(z_0+r \, e^{i
  t}\right) \, dt $$
for all $r>0$ small enough. Hence $w_2$ is subharmonic on $K_R \backslash \{ 0
\}$.

\smallskip

Further, condition (iv) implies 
$$\lim \limits_{r \to 0}
  \frac{M_{w_1}(r)}{\log(1/r)}=\lim \limits_{r \to 0}
  \frac{M_{u_2}(r)}{\log(1/r)}<+\infty$$

and  $w_2$ has therefore a subharmonic extension to $K_R$ by Lemma
\ref{lem:sub1}.  Thanks to the boundary behaviour of
 $u_1$ and $u_2$, {\it cf}.~assumption (iii), we have
$$\limsup_{z\to \xi} w_2(z) \le 0 \quad \text{for every}\, \, \,   \xi \in  \partial K_R\,.$$
Applying the maximum principle for subharmonic functions to $w_2$, we deduce that $w_2 \le 0$ in $K_R$ and so $u_1 \le u_2 $ in $K_R \backslash \{ 0\}$.\hfill{$\blacksquare$}

\bigskip

The following two examples illustrate that neither the assumption
``$u_2$ is  subharmonic'' in condition (i) nor the assumption 
$\lim_{r \to 0}
  M_{u_2}(r)/\log(1/r)<+\infty$  in condition (iv) of
 Theorem \ref{thm:ext_max} can be dropped.

\begin{example}
 We pick $ \kappa \equiv 0$, and choose 
$u_1 \equiv 0$ and 
$$ u_2: \overline{\D}\backslash\{ 0 \} \to \R\, , \quad\,  z \mapsto -\left( \frac{ \Re(z)}{|z|}+1 \right) \frac{1}{|z|^{3/2}}
  \log \frac{1}{|z|}\, . $$
It is easy to check that
the function $u_2$ is superharmonic in $\D \backslash \{ 0 \}$ and therefore a 
supersolution to $\Delta u=-\kappa(z) \, e^{2 u}$. 
Obviously, $u_1$ is a subsolution to $\Delta u=-\kappa(z) \, e^{2 u}$ in $\D$ and $u_1\le u_2$ on
$\partial \D$. However $u_1 \not \le  u_2$  in $\D\backslash \{ 0 \}$, although  
  $$ \lim \limits_{r \to 0}
  \frac{M_{u_1}(r)}{\log(1/r)}=\lim \limits_{r \to 0}
  \frac{M_{u_2}(r)}{\log(1/r)}=0 \, .$$
\end{example}

\begin{example}
Here we set $\kappa \equiv -e^2$ and consider on $\D \backslash \{ 0 \}$ the harmonic function 
$u_2(z)=\Re(z)/|z|^2$,
where
$$\inf_{|z|=1} u_2(z)=-1 \quad \text{and} \quad \lim_{r \to 0}\frac{M_{u_2}(r)}{\log(1/r)}=+\infty\, .$$
For the function $u_1: \D \backslash \{ 0 \} \to \R$ we choose 
$$u_1(z)= \log\left(\frac{1}{e} \frac{1}{|z| \log(e/|z|)}  \right)\, .$$
Thus  $u_1$ is a solution to the ``boundary value problem''
$$\Delta u=e^2\, e^{2u} \, \,  \text{ in }  \D \backslash\{ 0 \} \,, \quad
u\equiv -1 \, \text{ on } \partial \D$$
which satisfies
$$\lim_{r \to 0} \frac{M_{u_1}(r)}{ \log(1/r)}=1\, . $$ 
Clearly, $u_1\not \le
u_2$  in $\D \backslash \{ 0 \}$.
\end{example}

\section{Main results and proofs}\label{sec:results}

We are now prepared to prove Theorem \ref{thm:hoelder}. We shall see 
in \S \ref{par:conti} that
the representation formulas (\ref{eq:hoelder}) and (\ref{eq:hoelder2})
as well as the continuity properties of the remainder functions $v$ and $w$
follow  from the extended maximum principle (Theorem
\ref{thm:ext_max}) and the potential theory of Section \ref{par:potential}
by constructing suitable sub-- and supersolutions with {\it
  constant} curvature. In order to describe the precise behaviour of the 
 derivatives of the remainder functions, however, one needs to find appropriate
sub-- and supersolutions with {\it variable} curvature; {\it cf}.~\S
\ref{par:abl1}--\S \ref{par:abl2}.
\medskip

It is convenient to
introduce the following notion.

\begin{definition}
A real--valued function $\kappa$ on a set $G \subseteq \C$ is called strictly
negative on $G$ if $-a \le \kappa(z) \le -A$ in $G$ for finite
constants $a>0$ and $A>0$. A function $\kappa : \D \backslash \{ 0 \}
\to \R$ is called strictly negative at $z=0$, if $\kappa$ is strictly negative
in some punctured neighborhood of $z=0$.
\end{definition}

\subsection{Classification of the isolated singularities and continuity of the
remainder functions} \label{par:conti}

We start with a discussion of the behaviour of the 
solutions to $\Delta u=-\kappa(z) \, e^{2  u}$ under the assumption
that the curvature function $\kappa$ is merely strictly negative at $z=0$.

\begin{theorem}\label{thm:beschraenkt}
Let $\kappa : \D \backslash \{ 0 \} \to \R$ be strictly negative at $z=0$ and 
 $u: \D \backslash \{ 0 \}
\to \R$  a $C^2$--solution to
$\Delta u = - \kappa (z) \, e^{2u}$ in $\D \backslash \{ 0
\}$. Then $\alpha \le 1$, where $\alpha$ is the order of $u$, i.e.
\begin{equation} \label{eq:order}
\alpha:=\lim \limits_{r \to 0} \frac{M_u(r)}{\log (1/r)} 
\end{equation}
and
\begin{alignat}{3}
u(z)&= -\alpha \, \log{|z|} + v(z)\, , \quad \quad \qquad & \text{ if } \alpha <1\, , \label{eq:repr1}\\[2mm]
u(z)&= - \log|z|-\log{\log{(1/|z|)}}+w(z)\, ,\quad \quad \qquad & \text{ if } \alpha=1\, ,  \label{eq:repr2}
\end{alignat}
where $v(z)$ is continuous at $z=0$  and $w(z)=O(1)$ as $z \to 0$.
\end{theorem}

Note that the limit $\alpha$ in (\ref{eq:order}) always exists and $\alpha>-\infty$
since $u$ is subharmonic 
in a punctured neighborhood of $z=0$ (see (\ref{eq:sub1})). 

\begin{remark} \label{rem:mcowen}
Theorem \ref{thm:beschraenkt} was proved before by Heins
  \cite{Hei62} for $\kappa(z) \equiv -4$. It was observed by
 McOwen \cite{McO93} that Heins' method can also be used to prove  Theorem
 \ref{thm:beschraenkt} in the general case.
We give a different proof below which is easier in many respects.
 This proof also illustrates in a simple situation
the technique we use to establish Theorem \ref{thm:hoelder}. Moreover, it
allows us to say more about the remainder function $w(z)$ in (\ref{eq:repr2})
under some mild extra
assumptions on the curvature function $\kappa(z)$; {\it cf}.~Theorem
\ref{thm:alpha=1,continuity}. These extra information will be crucial in
geometric applications such as Theorem \ref{thm:2} and Theorem \ref{thm:3}.
 An inspection of the method
of Heins and McOwen shows that their technique does not seem to be  capable of
yielding  such
refinements of Theorem \ref{thm:beschraenkt}. 
We also note that Theorem \ref{thm:beschraenkt} only
deals with the case of strictly negative curvature. Positive
curvature functions appear to be more difficult to handle
 and some partial results in this case like one--sided estimates 
have recently been obtained by Yunyan \cite{Yun03}.
\end{remark}

{\bf Proof of Theorem \ref{thm:beschraenkt}.}

Without loss of generality we may assume $u \in C(\overline{\D} \backslash \{
0 \})$ and 
\begin{equation} \label{eq:b}
-a \le \kappa(z) \le -A \qquad \text{ in } \D \backslash \{ 0 \}
\end{equation}
for some finite constants $a>0$ and $A>0$. We first show that
the order  $\alpha$ of $u(z)$ at $z=0$ is always $\le 1$. To see this just
note that $u$ is a subsolution to $\Delta \nu=A \, e^{2 \nu}$ in $\D
\backslash \{ 0 \}$, so
\begin{equation}\label{eq:max_supersol}
 u(z)\le \log\left( \frac{1}{\sqrt{A} }\frac{1}{|z| \log(1/|z|)}  \right)\,,
 \qquad z \in \D \backslash \{ 0 \} \, , 
\end{equation}
because the function on the right--hand side is the maximal solution
to $\Delta \nu=A \, e^{2 \nu}$ in $\D \backslash \{ 0 \}$. This  follows from Ahlfors' lemma \cite{Ahl38}
by noting that 
$$ \frac{|dz|}{\sqrt{A} \, |z| \, \log(1/|z|)}$$
is the hyperbolic metric of the punctured disk $\D \backslash \{ 0 \}$
with constant curvature $-A$. Now, (\ref{eq:order}) is an immediate consequence of
the inequality (\ref{eq:max_supersol}).

\medskip

In a next step we construct a supersolution to $\Delta u=-\kappa(z) \, e^{2
  u}$ by looking for a solution $u^A_{\alpha}$ 
to $\Delta \nu=A \, e^{2 \nu}$ with order $\alpha$  at $z=0$. Indeed it is not difficult to show that

\begin{equation} \label{eq:help1}
u_{\alpha}^A(z):=
\begin{cases}  \displaystyle\log \left( \frac{2}{\sqrt{A}} \frac{(1-\alpha) \,
    |z|^{-\alpha}}{1-|z|^{2(1-\alpha)}} \right) & \hspace{0.5cm} \text{if} \,
        \,  \,  \alpha<1\,, \\[6mm]
        \displaystyle  \log \left( \frac{1}{\sqrt{A}} \, \frac{1}{  |z| \,
        \log (1/|z|) }\right) & \hspace{0.5cm}\text{if} \, \, \, 
  \alpha=1\, , 
\end{cases} 
\end{equation}
has the required properties. For $\alpha <1$ 
this supersolution is obtained by noting that 
the conformal Riemannian metric
$$\lambda_{\alpha}^A(z)\, |dz|:= e^{u_{\alpha}^A(z)}\, |dz|$$ 
is
 the formal pullback of the hyperbolic metric of the unit disk with constant curvature
  $-A$
$$\frac{2}{\sqrt{A}} \,\frac{|dz|}{1-|z|^2}\, ,$$
under the map $z \mapsto z^{1-\alpha}$. We are now in a position to apply the
extended maximum principle (Theorem \ref{thm:ext_max}) which leads to
\begin{equation}\label{eq:super}
u(z)\le u_{\alpha}^A(z) \quad \text{for} \, \, \, z \in \D \backslash \{ 0
\}\, .
\end{equation}
In order to get a lower bound for $u(z)$, we next look for appropriate
subsolutions to $\Delta u=-\kappa(z) \, e^{2
  u}$, i.e., in view of (\ref{eq:b}), for solutions 
to $\Delta \nu=a \, e^{2 \nu}$  with order $\alpha$ at $z=0$.
A one--parameter family of such solutions $u^{a}_{\alpha,R}$\,, $R>1$,
is obtained from the functions $u^A_{\alpha}$ in (\ref{eq:help1}) by
replacing $A$ with $a$ and a suitable rescaling:
$$ u^a_{\alpha,R}(z):=u^a_{\alpha}(z/R)+\log (1/R) \, .$$
Geometrically, the conformal Riemannian metric $e^{u^a_{\alpha,R}(z)} \, |dz|$ is the pullback
of the hyperbolic metric
\begin{itemize}
\item[(i)] 
of the unit disk $\D$ with constant curvature $-a$ under the map $z \mapsto (z/R)^{1-\alpha}$  if
$\alpha<1$;
\end{itemize}
and
\begin{itemize}
\item[(ii)]
of the punctured unit disk $\D \backslash \{ 0 \}$ with constant curvature
$-a$ under the map $z \mapsto z/R$ if $\alpha=1$.
\end{itemize}
Thus, $u^a_{\alpha,R}$ is a subsolution to $\Delta u=-\kappa(z) \, e^{2 u}$ in
$\D \backslash \{ 0 \}$ with order $\alpha$ at $z=0$ for each $R>1$.
We now have to choose the parameter $R>1$ in an appropriate way. For this denote
$$ m:=\min \limits_{|z|=1} u(z) \, . $$
Observe that for every $|z|=1$
$$ \lim \limits_{R \to \infty} u^a_{\alpha,R}(z)=
\lim \limits_{R \to \infty} u^a_{\alpha,R}(1)=-\infty \, .$$
Consequently, we can find $R>1$ such that $u^a_{\alpha,R}(z) \le m$ for $|z|=1$.
Hence, we can again apply the extended maximum principle and obtain
\begin{equation}\label{eq:sub}
u_{\alpha,R}^a(z) \le u(z) \quad \text{for} \, \, \, z \in \D \backslash \{ 0
\}\, .
\end{equation}
Combining (\ref{eq:super}) and (\ref{eq:sub})  yields 
$$
\displaystyle \log \left(\frac{2(1-\alpha )}{\sqrt{a}} \,
\frac{1}{R^{1-\alpha}(1-|z/R|^{2(1-\alpha)})} \right) \le \, u(z)+ \alpha \log|z|\, 
\le \log\left( \frac{2(1-\alpha )}{\sqrt{A}} \, \frac{1}{1-|z|^{2(1-\alpha)}}
\right)$$ for $\alpha <1$
and 
$$\log \left( \frac{1}{\sqrt{a}} \,
\frac{\log{(1/|z|)}}{\log{(R/|z|)}} \right) \le \,
u(z) + \log|z| + \log\log \frac{1}{|z|}\, 
\le \log \left(\frac{1}{\sqrt{A}} \right)$$
for $\alpha=1$. Thus $u$ has the desired representation (\ref{eq:repr1}) and
(\ref{eq:repr2}) respectively, where the remainder functions $v$ and $w$ are
continuous in $\D \backslash \{ 0 \}$ and
bounded at $z=0$.

\medskip

To complete the proof of Theorem \ref{thm:beschraenkt}, we need to show that
for $\alpha<1$ the remainder function $v(z)$ is in fact continuous at $z=0$.
Using the PDE $\Delta u=-\kappa(z) \, e^{2 u}$ we get that $v(z)=u(z)+\alpha
\log |z|$ is a solution to
$$\Delta v= \frac{- \kappa(z)}{|z|^{2\alpha}} \, e^{2 v}  \qquad \text{for} \,
\, z \in \D \backslash \{ 0 \}\, ,$$
which shows
that $v$ is subharmonic on $\D \backslash \{ 0 \}$. On the other hand $v$ 
 is bounded  in a
neighborhood of $z=0$ and consequently subharmonic in all of $\D$. The
fact that $z \mapsto \Delta v(z)$ is integrable over $K_r$ for each $0<r<1$
 legitimizes the use of
Proposition \ref{thm:poisson-jensen} and leads to
\begin{equation}\label{eq:representation1}
v(z)=h(z)+\frac{1}{2\pi} \iint \limits_{K_r} \log|z- \xi| \,
\frac{-\kappa(\xi)}{|\xi|^{2\alpha}}\,  e^{2v(\xi)} \, d\sigma_{\xi}\, ,
\quad z \in K_r\, ,
\end{equation}
where $h$ is a harmonic function on $K_r$.
The continuity of $v$ in all of $\D$ is now an immediate consequence of Proposition \ref{thm:potential}.\hfill{$\blacksquare$}

\bigskip

In Theorem \ref{thm:beschraenkt} only the fact that $\kappa$ is
strictly negative at $z=0$ 
guarantees the continuity of the remainder function $v$ at $z=0$ for $\alpha<1$. 
This is no longer true if $\alpha=1$ as the following example shows.

 \begin{example}\label{ex:alpha=1,kappa=beschraenkt}
The function 
$$u(z)= -\log|z| - \log \log \frac{1}{|z|}+2+ \frac{\sin \left(
    \log\log(1/|z|)\right)}{6+\sin \left(
    \log\log(1/|z|)\right) }
$$
is a solution to
$$\Delta u = - \kappa(z)\, e^{2u}$$
in $\D \backslash  \{ 0 \}$, where
$$\kappa(z)= \left( 6 \, \frac{\sin \beta(z)+\cos\beta(z)}{(6+\sin
    \beta(z))^2}+12 \, \frac{(\cos \beta(z))^2}{(6+\sin
    \beta(z))^3}-1 \right) \,  \exp \left(-6+\frac{12}{6+\sin(\beta(z))}
\right) \, ,$$
with $\beta(z)=\log\log(1/|z|)$.
Obviously, both the curvature $\kappa(z)$ and the remainder function 
$w(z)$ are bounded
when $z \to 0$ but not continuous at
$z=0$. 
\end{example}

Thus the case $\alpha=1$ is exceptional as far as the continuity of the
remainder function $w$ is concerned.
We now show that  the remainder function
$w$ is continuous at $z=0$, when we additionally assume  that the curvature $\kappa$
is continuous at $z=0$.

\begin{theorem}\label{thm:alpha=1,continuity}
Let $\kappa :\D \to \R$ be a continuous function with $\kappa(0)
<0$. If $u(z)= -\log|z| -\log \log (1/|z|) +w(z)$, where $w(z)=O(1)$ for $z
\to 0$, is  a solution to $\Delta u=-
\kappa(z)\, e^{2u}$ in $\D \backslash \{ 0 \}$, then $w$ is continuous at $z=0$ and $w(0)=-\log
\sqrt{-\kappa(0) }$.
\end{theorem}

We wish to point out that in contrast to the cases $\alpha<1$, where the
continuity of the remainder function $v$ follows from potential--theoretical
considerations under weaker assumptions, this method does not work in the
exceptional case $\alpha=1$,
even though there is  a corresponding representation formula for the remainder
function $w$  (see formula (\ref{eq:representation2}) below).
Instead, we use the method of super-- and subsolutions and the generalized
maximum principle  (Theorem \ref{thm:ext_max}) 
for this purpose.

\medskip

{\bf Proof of Theorem \ref{thm:alpha=1,continuity}.}\\
Set $b:= -\kappa(0)$ and choose $0< \varepsilon < b/2$. By the
continuity of $\kappa$ there is a disk $K_{\varrho}$ such that
$$0 < b- \varepsilon \le - \kappa(z) \le b+ \varepsilon \qquad \text
{for} \, \, z \in K_{\varrho}\, .$$
Next, define
$$R=\varrho \, \exp\left(1/\left(\varrho\, \sqrt{b}\,
    \tilde{m}_{\varrho}\right)\right)
 \quad \text{and} \quad
R'=\varrho \, \exp\left(1/\left(\varrho\, \sqrt{b}\, \tilde{M}_{\varrho}\right)\right)\, ,$$
where
$$\tilde{M}_{\varrho}= \max_{|z|=\varrho} e^{u(z)} \quad \text{and}  \quad \tilde{m}_{\varrho}= \min_{|z| =
  \varrho} e^{u(z)}\, ,$$
and let
$$ u_{-}(z):=\log \left( \frac{1}{\sqrt{b+\varepsilon}}\,  \frac{1}{|z| \, \log
     (R/|z|)}\right) \qquad \text{ and } \qquad u_+(z):=\log\left(
   \frac{1}{\sqrt{b-\varepsilon}}\,  \frac{1}{|z| \, \log( R'/|z|)} \right) \, .$$
Then $u_{-}$ is a subsolution to $\Delta u=-\kappa(z) \, e^{2 u}$ in $K_{\varrho}
\backslash \{ 0 \}$
with order 
$1$ at $z=0$ and $u_{-} \le u$ on $|z|=\varrho$. Also,
$u_+$ is a supersolution to $\Delta u=-\kappa(z) \, e^{2 u}$ in $K_{\varrho} \backslash
\{ 0 \}$
with order $1$ at $z=0$ and $u \le u_+$ on $|z|=\varrho$.
By the extended maximum principle (Theorem \ref{thm:ext_max}), we get
$$\log \left( \frac{1}{\sqrt{b+\varepsilon}}\,  \frac{1}{|z| \, \log (R/|z|)}\right) \le u(z) \le
\log\left(  \frac{1}{\sqrt{b-\varepsilon}}\,  \frac{1}{|z| \, \log( R'/|z|)} \right) \quad \text{for} \,
    z \in K_{\varrho} \backslash\{ 0 \}\, .$$
Rearranging the latter inequality yields
$$w(z) \to \log \left( \frac{1}{\sqrt{b}} \right) \qquad \text{for} \, \, z \to 0\, ,$$
which is what had to be proved.
\hfill{$\blacksquare$}

\subsection{First derivatives of the remainder functions: Results}  \label{par:abl1}

We now turn to a discussion of the properties of the first partial derivatives
of the remainder functions $v$ and $w$. We first consider the case of strictly
negative curvature functions.

\begin{theorem}\label{thm:firstderivative}
Let $\kappa: \D\backslash \{ 0 \} \to \R$ be strictly negative at $z=0$
and let $u: \D\backslash \{ 0 \} \to \R$
be a $C^2$--solution to 
$\Delta u =-\kappa(z) \, e^{2u}$
in $\D \backslash \{ 0 \}$ with order $\alpha \le 1$ at $z=0$, i.e.,
 $u(z):= -\alpha \log|z| +v(z)$ for $\alpha <1$ and
$u(z):=-\log|z|-\log\log (1/|z|)+w(z)$ for $\alpha =1$.  Then
\begin{alignat*}{2}
&v_{z}(z),\, v_{\overline{z}}(z)\, \,  \text{are continuous at}
\, \,  z=0\qquad \quad &\text{if}\, \, \,  &\alpha < 1/2;\\[-6mm]
\intertext{and}\\[-10mm]
&v_{z}(z),\, v_{\overline{z}}(z)=O(\log(1/|z|))\qquad  &\text{if}\, \, \, &\alpha =1/2,\\[1mm]
&v_{z}(z),\, v_{\overline{z}}(z)=O\left(|z|^{1-2\alpha}\right) \qquad &\text{if}\, 
\, \, &1/2<\alpha <1,\\[1mm]
&w_{z}(z),\, w_{\overline{z}}(z)=O\left(|z|^{-1} (\log(1/|z|))^{-1}\right) \qquad & \text{if}\,
\, \,  &\alpha =1,
\end{alignat*}
when $z$ approaches $z=0$.
\end{theorem}

\begin{remark}
\begin{itemize}
\item[(a)]
Theorem \ref{thm:firstderivative} is sharp. In fact, for 
$\alpha\not \in \{ 1/2 ,1\}$, this is illustrated with Example \ref{ex:dani1}.
For $\alpha=1/2$ see Example \ref{ex:3.6}, and for $\alpha=1$ one can use the
function of Example \ref{ex:alpha=1,kappa=beschraenkt}.
\item[(b)]
We also note that  if $\alpha \not=1/2$ and $\alpha\not=1$, then Theorem
  \ref{thm:firstderivative} refines (part of) Theorem \ref{thm:hoelder}.
\item[(c)] In view of (a) and (b) we see that
 $\alpha=1/2$
  and $\alpha=1$ are exceptional cases as far as the first derivatives of the
  remainder functions $v$ and $w$ are concerned. Also,
 in order to prove 
the (stronger) assertions of Theorem
  \ref{thm:hoelder}  for the first derivatives of $v$ and $w$ in those
  exceptional cases, it does not suffice 
to assume that $\kappa$ is  strictly negative at $z=0$. As we shall see
in Example \ref{ex:3.7} (for $\alpha=1/2$) and Example \ref{ex:3.10} (for
$\alpha=1$) it is  even not enough to suppose that $\kappa$ is continuous (and
negative) at $z=0$. 
\end{itemize}
\end{remark}

In the following theorem we therefore adopt the 
standard assumption for second order elliptic PDEs and consider
curvature functions $\kappa$ which are  H\"older continuous on $\D$.
This assumption  was
used before for the curvature equation for instance by Troyanov 
\cite[p.~800]{Tro90}, but see also the remarks following the next theorem.

\begin{theorem} \label{thm:exceptional}
If $\kappa$ is locally H\"older continuous in $\D$ with $\kappa(0)<0$, then
with the notation of Theorem \ref{thm:firstderivative},
\begin{alignat*}{2}
&v_{z}(z),\, v_{\overline{z}}(z)=O(1)\quad  &\text{if}\, \, &\alpha
=1/2\, ,\\[-6mm]
\intertext{and}\\[-12mm]
&w_{z}(z),\, w_{\overline{z}}(z)=O\left(|z|^{-1}\,( \log(1/|z|))^{-2}\right) \quad \, \, &\text{if}\,
\, &\alpha =1\,,
\end{alignat*}
when $z$ tends to $z=0$.
\end{theorem}

By Example \ref{ex:dani1}, Theorem \ref{thm:exceptional} is best possible.
We note that our proof of the case $\alpha=1/2$ in Theorem \ref{thm:exceptional}
only uses  potential--theoretic tools. For the case $\alpha=1$, by contrast,
potential theory alone is not enough. We require further {\it a--priori}
information about the solutions to the curvature equation in this case and 
we are going to use our extended maximum principle for this purpose.
However, in order to be able to
apply the extended maximum principle we need to construct suitable sub-- and
supersolutions to $\Delta u=-\kappa(z) \, e^{2 u}$, which take 
the H\"older continuity of the curvature function $\kappa$ into account. This is the key
point. It turns out that it is easier to construct the required sub--
and supersolution under weaker assumptions on the curvature
function.\footnote{Die Vergr\"o{\ss}erung der Beweislast kann also vorteilhaft
  sein: denn sie st\"arkt den
  Beweistr\"ager \cite{Pol}.} We have the following crucial lemma.

\begin{lemma}\label{lem:wachstum}
Let $\kappa: \D \to \R$ be a function with $\kappa(0)<0$ and
\begin{equation}\label{eq:kruemmung_vergleich}
\kappa(z)=\kappa(0) + \frac{r(z)}{ (\log(1/|z|))^{2}}\, ,
\end{equation}
where $r(z)=O(1)$ as $z \to 0$. If $u: \D \backslash \{ 0 \} \to \R$ is a solution to
$\Delta u= - \kappa(z) \, e^{2u}$
in $\D \backslash \{ 0 \}$ with 
$u(z)= - \log|z|-\log\log(1/|z|) +w(z) $
where $w(z)=O(1)$ for $z\to 0$,
then  
there is a disk $K_{\varrho}$ such that
\begin{equation}\label{eq:q}
 \displaystyle \left| -\kappa(z)\,  e^{2w(z)}-1 \right| \le
\frac{C}{\log(1/|z|)} \, , \qquad z \in K_{\varrho} \, ,
\end{equation}
for some constant $C>0$.
\end{lemma}

Clearly, every H\"older continuous curvature function $\kappa : \D \to \R$
fulfills condition (\ref{eq:kruemmung_vergleich}) of Lemma \ref{lem:wachstum}.
The estimate (\ref{eq:q}) is exactly the {\it a--priori} information we need
to know about the remainder function $w$ in order to start up the
potential--theoretic
part of the proof of Theorem \ref{thm:exceptional} for $\alpha=1$.
We also note that Lemma \ref{lem:wachstum} is in a certain sense a higher--order
version of the classical Ahlfors' lemma \cite{Ahl38}.

\subsection{First derivatives of the remainder functions: Proofs and examples}
\label{par:proofs}

We now give the proofs of the results of Paragraph \ref{par:abl1} and illustrate
 their sharpness with a number of examples. It suffices to
restrict ourselves to the investigation of the
derivative with respect to $z$.

\medskip

We start with the proof of Theorem \ref{thm:firstderivative} for $\alpha<1$.

\medskip

{\bf Proof of Theorem \ref{thm:firstderivative} for $\boldsymbol{\alpha<} \mathbf{1}$.}

Without loss of generality we may assume $-a  \le \kappa(z)\le 
 -A$ in $\D$, where $a$ and $A$ are positive constants. We now make essential use
 of the representation formula (\ref{eq:representation1}) for the remainder
 function $v(z)$, which has been established in the course of the proof
of Theorem \ref{thm:beschraenkt}.

\smallskip

If  $\alpha <1/2$  the claim is an immediate consequence of  (\ref{eq:representation1}) and Proposition \ref{thm:potential}.

\smallskip

For $1/2\le \alpha < 1$, pick $R<1$. Then formula (\ref{eq:representation1})
combined with Proposition \ref{thm:potential} gives
$$|v_z(z)| \le \sup_{z \in K_R} |h_z(z)|+ \sup_{z \in K_R}\left(
  -\kappa(z)\, 
  e^{2v(z)}\right) \, \frac{1}{2\pi} \iint \limits_{K_R} \frac{1}{2}\frac{1}{|z-\xi|}\, 
  \frac{1}{|\xi|^{2\alpha}} \, d\sigma_{\xi}\, , \qquad z \in K_R \backslash
  \{0 \}\, .$$
where $h$ is a harmonic function on $K_R$.
To see that
$$
v_z(z)=
\begin{cases}
O\left(\log(1/|z|)\right) \quad &\text{if} \, \, \, \alpha=1/2, \\[2mm]
O\left(|z|^{1-2\alpha}\right) \quad &\text{if} \, \, \, 1/2< \alpha<1,
\end{cases}
$$
we just note that 
$$\frac{1}{2\pi} \iint \limits_{K_R} \frac{1}{|z-\xi|}\, 
  \frac{1}{|\xi|} \, d\sigma_{\xi} \le C_1 +\log \frac{1}{|z|}\,, \qquad
 \text{if} \, \, \,  \alpha=\frac{1}{2}\, ,$$
 and
$$\frac{1}{2\pi} \iint \limits_{K_R} \frac{1}{|z-\xi|}\, 
  \frac{1}{|\xi|^{2\alpha}} \, d\sigma_{\xi}\le C_2 \, |z|^{1-2\alpha}\,, \qquad
  \text{if} \, \, \, 
  \frac{1}{2}<\alpha <1\, ,$$
(see \cite[p.~215]{Els91}) for some positive constants $C_1$ and
$C_2$. \hfill{$\blacksquare$}

\bigskip

If the curvature $\kappa$ is strictly negative at $z=0$, but not necessarily
continuous, then in the preceding proof
the estimate for the case $\alpha=1/2$  is best possible:

\begin{example} \label{ex:3.6}
Let $\kappa: \D \backslash \{ 0 \} \to [-6\, e,-4\, e^{-8}]$ be given by
$$ z \mapsto -\left( 4 + 2\,  \frac{|\Re
  (z)|}{ |z|} \right)\, \exp\left(-2 \, \left( |\Re(z)| \, \log|z| +4\, |z|  \right)
  \right)\, .$$
 Then the function $u(z) = -1/2 \log |z|+v(z)$ where $v(z)= |\Re(z)|\, \log
  |z|+ 4 |z|$ is a solution to
$$\Delta u= -\kappa(z)\,  e^{2u} \qquad \text{in} \, \, \D \backslash \{ 0 \}$$
and
$$v_z(z)= \frac{|\Re(z)|}{2 z} +2 \, \frac{\overline{z}}{|z|}- \frac{
  \Re(z)}{2\, |\Re(z)|} \log(1/|z|)\, .$$
\end{example}

\medskip

{\bf Proof of Theorem \ref{thm:exceptional} for $\boldsymbol{\alpha=}\mathbf{1/2}$.}

 At first we observe the assumption that $\kappa$ is locally H\"older
 continuous in $\D$ (without loss of generality with fixed exponent $\gamma$,
 $0< \gamma \le 1$, say)
 implies that $v$ is locally H\"older continuous in $\D$ . To this end 
 we may suppose $-a \le\kappa(z)\le-A$ in $\D$, where $a$ and $A$ are
positive constants. Further, we  define in
$\D \backslash \{ 0 \}$ the regular conformal Riemannian  metric $\mu(z)\, |dz|=e^{u(z)}\, |dz|$ and
write $\mu(z)=|z|^{-1/2} \sigma(z)$, so
 $\sigma(z)=e^{v(z)}$  is continuous at $z=0$ and $\sigma (0)>0$;
see Theorem \ref{thm:beschraenkt}. Now
$$\tilde{\mu}(z)\, |dz|:= 2 \, \mu(z^2)\, |z| \, |dz|=2 \, \sigma(z^2) \, |dz|$$
is a (positive) continuous conformal Riemannian metric in $\D$ which is regular $\D\backslash \{ 0 \}$. Its curvature
$\kappa_{\tilde{\mu}}$ is locally H\"older continuous in $\D$
 as
$\kappa_{\tilde{\mu}}(z)=\kappa(z^2)$.
Thus $\tilde{\mu}(z)$ has a $C^2$--extension to $\D$
by elliptic regularity.
In particular, $\tilde{\mu}_z$ is bounded in $K_R$ for some $R<1$.
By construction, this implies
$$|z \cdot v_z(z)| \le C\, |z|^{1/2}$$ 
for $z \in K_R$  and some constant $C$, so $v$
 is locally H\"older continuous in $\D$ with exponent $1/2$.

\smallskip

This in turn implies that the function $- \kappa(\xi)\, e^{2v(\xi)}$ is locally H\"older continuous in $\D$ with
exponent $\tilde{\gamma}=\min\{\gamma,1/2\}$ and 
 $$\kappa(\xi)\, e^{2v(\xi)}= \kappa(0) \, e^{2v(0)} +r(\xi)\, |\xi|^{\tilde{\gamma}}\quad \text{ for} \, \,  \xi \in \D\, , $$
 where $r$ is a continuous function on $\D\backslash \{ 0  \}$ with
 $r(\xi)=O(1)$ for $\xi \to 0$.
Thus we have for $z \in K_R \backslash \{ 0 \}$ in view of formula
(\ref{eq:representation1}) and  Proposition \ref{thm:potential}
\begin{equation}\label{eq:alpha=1/2,kappa=hoelder_1}
v_z(z)= h_z(z)+ \frac{1}{2\pi} \iint\limits_{K_R} \frac{1}{2(z-\xi)}
\frac{ - \kappa(0) \, e^{2v(0)}}{|\xi|} \, d\sigma_{\xi}- \frac{1}{2\pi} \iint\limits_{K_R} \frac{1}{2(z-\xi)}
\frac{r(\xi)}{|\xi|^{1- \tilde{\gamma}}}\, d\sigma_{\xi}\, ,
\end{equation} 
where $h$ is harmonic on $K_R$.

\smallskip

To confirm that $v_z(z)=O(1)$ as $z \to 0$, we first observe that the second integral of
(\ref{eq:alpha=1/2,kappa=hoelder_1}) is bounded at $z=0$. 
For the first integral of (\ref{eq:alpha=1/2,kappa=hoelder_1}) we define $k:= -
\kappa(0)\, e^{2v(0)}$ and consider  the subharmonic function $\omega(z)= k
|z|$ on $\D$. Since $\Delta \omega(z)=k/ |z|$ for $z \in \D \backslash \{ 0 \}$
we have by Proposition \ref{thm:poisson-jensen}
$$\omega(z)= h_{\omega}(z)+\frac{1}{2\pi} \iint \limits_{K_R} \log|z-\xi|\, 
\frac{k}{ |\xi|}\, d\sigma_{\xi}\, ,\qquad z \in K_R\, ,$$
where $h_{\omega}$ is harmonic on $K_R$.  Lastly, differentiating yields
$$\omega_z(z)=k \frac{\overline z}{|z|} \quad \text{and} \quad \omega_z(z)=  h_{{\omega}_z}(z)+\frac{1}{2\pi} \iint
\limits_{K_R} \frac{1}{2\, (z - \xi)} \frac{k}{|\xi|} \, d\sigma_{\xi}\, $$
for $z \in K_R\backslash \{ 0 \}$. A comparison with
(\ref{eq:alpha=1/2,kappa=hoelder_1})  completes the proof. \hfill{$\blacksquare$}

\bigskip

For $\alpha=1/2$
we note that if $\kappa$ is merely continuous at $z=0$, then $v_z$ is not
necessarily bounded at $z=0$:

\begin{example} \label{ex:3.7}
Consider the function $u: \D \backslash \{ 0 \} \to \R$ defined by
$$ u(z)=-\frac{1}{2} \, \log  |z| +\Re(z)  \, \log \log(e/|z|)+ 4\, |z|\, .$$ 
Then $u$ is a solution to
$$\Delta u=-\kappa(z)\, e^{2u}$$
in $\D \backslash \{ 0 \}$, where
$$\kappa(z)=-\left( 4+ \frac{\Re(z)}{|z|}\, \frac{-3+2 \log|z|}{(1-
    \log|z|)^2} \right)\,  \exp\left( -2 \left( \Re(z) \log
    \log(e/|z|) +4\, |z| \right)  \right) \, .$$
Note that $\kappa(z)$ 
is continuous in $\D$ with $\kappa(0)=-4$ but not locally H\"older continuous at $z=0$.
For the derivative of 
$$v(z)= \Re(z)  \, \log \log(e/|z|) + 4\, |z|$$
we obtain
$$v_z(z)=2\frac{\overline{z}}{|z|} - \frac{1}{2} \frac{\Re(z) \, 
    \overline{z}}{|z|^2 \log(e/|z|)}+ \frac{1}{2} \log \log (e/|z|) \, .$$
Thus $|v_z(z)| \to \infty$ for $z \to 0$.
\end{example}

{\bf Proof of Theorem \ref{thm:firstderivative} for $\boldsymbol{\alpha=1}$.}

(i) \,  We first show that for every $0<r<1$
\begin{equation}\label{eq:representation2}
w(z)=h(z)+\frac{1}{2\pi} \iint \limits_{K_r} \log|z- \xi| \,
\frac{-\kappa(\xi)  e^{2w(\xi)} -1}{|\xi|^{2}\,  (\log(1/|\xi|))^2}\,  \,
d\sigma_{\xi} \qquad \text{for} \, \, z \in K_r\backslash \{ 0 \} \, ,
\end{equation}
 where $h$ is a harmonic function on $K_r$.
This is the analogue  to the corresponding
representation formula (\ref{eq:representation1}) for the function $v$ for the
cases $\alpha<1$.

\smallskip 

In order to prove (\ref{eq:representation2})
we define for $z \in K_r \backslash \{ 0 \}$ the function
$$p(z):= w(z)-\log \log (1/|z|) \, .$$
Since
$$\Delta p(z)= \frac{- \kappa (z)}{|z|^2 \, (\log(1/|z|))^2} \,
  e^{2w(z)} \, >0 \quad \text{for} \, z \in K_r \backslash \{ 0 \} $$
and $\lim_{z \to 0} p(z) = -\infty$, $p$ is subharmonic on $K_r$. 

\smallskip

 By Proposition \ref{thm:poisson-jensen}, as $z \mapsto \Delta p(z)$ is
  integrable over $K_r$, we obtain
$$p(z)=h_p(z) +\frac{1}{2\pi} \iint \limits_{K_r} \log|z - \xi| \, \frac{
  -\kappa(\xi)\,  e^{2w(\xi)}}{|\xi|^2\, (\log(1/|\xi|))^2} \,\,
d\sigma_{\xi} \, ,\qquad 
z \in K_r \,, $$
 where $h_p$ is harmonic on $K_r$.

\smallskip 

Applying Proposition \ref{thm:poisson-jensen} again, this time to the subharmonic
function $q(z):=- \log \log (1/|z|)$, we deduce that 
$$q(z)=h_q(z) + \frac{1}{2\pi} \iint\limits_{K_r} \log|z- \xi| \, 
\frac{1}{|\xi|^2 \, (\log(1/|\xi|))^2} \, \, d\sigma_{\xi}\qquad \text{for} \, \,
z \in K_r\, ,$$ 
where $h_q$ is harmonic on $K_r$. This gives (\ref{eq:representation2}) with
the harmonic function $h(z)=h_p(z)-h_q(z)$.

\medskip

(ii) \, Let $R < 1/(2e^4)$ and set $q(z)= -\kappa(z)\, e^{2w(z)}
-1$. According to the representation formula (\ref{eq:representation2}) and
  Proposition \ref{thm:potential} we obtain
$$ |w_z(z)|\le \sup_{z \in K_R}|h_z(z)| + \frac{1}{2\pi} \iint \limits_{K_R}
  \frac{1}{2} \frac{1}{|z- \xi|}\, \,\frac{ |q(\xi)|}{|\xi|^{2}\, 
  (\log(1/|\xi|))^2}\, \, d\sigma_{\xi}\, , \qquad z \in K_R\backslash \{ 0
  \}\, ,
$$
where $h$ is a harmonic function on $K_R$.
 We only have to find that
for small $z$,
$$ I(z):=\frac{1}{2\pi} \iint \limits_{K_R} \frac{1}{|z- \xi|}\, \,\frac{|q(\xi )|}{|\xi|^{2}\, 
  (\log(1/|\xi|))^2}\, \, d\sigma_{\xi} \le \frac{C}{|z|\log(1/|z|)}
$$
for an appropriate constant $C$.
To this end fix $z \in K_{R/2}$ and set  $r=|z|/2$. Further,  let  $N= \{\xi \in K_R \backslash (K_r \cup K_r(z)): |z-\xi | \ge |\xi|
\, \}$
and $P= \{\xi \in K_R \backslash (K_r \cup K_r(z)):
|z-\xi | < |\xi|\,  \} $,
where $K_r(z)$ stands for the disk about $z$ with radius $r$. Then, if
$M=\sup_{\xi \in K_R} |q(\xi)|$,
\begin{align}
I(z)&\le
\frac{1}{2\pi} \iint \limits_{K_r} \frac{1}{|z- \xi|} \,\frac{|q(\xi)|}{|\xi|^{2}\, 
  (\log(1/|\xi|))^2}\, \, d\sigma_{\xi}
  + \frac{1}{2\pi} \iint \limits_{K_r(z)} \frac{1}{|z- \xi|} \,\frac{M}{|\xi|^{2}\,
 ( \log(1/|\xi|))^2}\, \, d\sigma_{\xi} \notag\\[2mm]
&\phantom{\le} \quad +\frac{1}{2\pi} \iint \limits_{N} \frac{1}{|z- \xi|}\, \,\frac{M}{|\xi|^{2}\, 
  (\log(1/|\xi|))^2}\, \, d\sigma_{\xi}
  +\frac{1}{2\pi} \iint \limits_{P} \frac{1}{|z- \xi|}\, \,\frac{M}{|\xi|^{2}\, 
  (\log(1/|\xi|))^2}\, \, d\sigma_{\xi} \notag\\[2mm]
&\le \sup_{\xi \in K_r}|q(\xi)|\, \, \frac{1}{r} \, \int_0^r \frac{1}{\varrho
 ( \log(1/\varrho))^2} d\varrho\,  + 
\frac{M}{r^2 \, (\log(1/r))^2} \int_0^r d\varrho\notag\\[2mm]
& \phantom{\le} \quad +\frac{M}{\sqrt{r} \, (\log(1/r))^2}
\int_r^R\frac{1}{\varrho^{3/2}} d\varrho  +\frac{M}{\sqrt{r} \, (\log(1/r))^2}
\int_r^{2R}\frac{1}{\varrho^{3/2}} d\varrho \notag\\[2mm]
& \le \sup_{\xi \in K_r}|q(\xi)|\, \, \frac{2}{|z|\, \log(1/|z|)}\, + \,
\frac{10M}{|z|\, (\log(1/|z|))^2}= \frac{C}{|z|\, \log(1/|z|)}
\label{eq:integral}
\end{align}
for some constant $C>0$.\hfill{$\blacksquare$}

\bigskip

For strictly negative $\kappa$ (not necessarily continuous)
 the above proof gives an optimal
estimate as it can be seen from Example \ref{ex:alpha=1,kappa=beschraenkt}.
If $\kappa$ is continuous in $\D$, then 
the function $q(z)=-\kappa(z) \, e^{2w(z)}-1$ is
continuous at $z=0$ with $q(0)=0$ by Theorem \ref{thm:alpha=1,continuity}. Thus $w_z(z)=o(|z|^{-1}
\log(1/|z|)^{-1})$ by (\ref{eq:integral}).
This is pretty sharp, because the following example shows that
 $|w_z (z)|=C\cdot  (|z|^{-1} (\log (1/|z|))^{-1-\beta})$ is possible for any $\beta>0$.

\begin{example} \label{ex:3.10}
Define for $\beta >0$ and  $z \in \D$ the continuous function
$$\kappa(z)=-\exp \left( \frac {-2}{ (\log(1/|z|))^{\beta} } \right) \left( 1
+\beta(1+\beta) \frac{1}{(\log(1/|z|))^{\beta} }\right)\, .$$
Then 
$$u(z)= - \log|z| - \log \log \frac{1}{|z|}  +
\frac{1}{\left(\log (1/|z|)\right)^{\beta}}\, , \quad  z \in \D \backslash \{ 0
\}\, ,$$ 
is a solution to
$$\Delta u = -\kappa(z) \, e^{2u}$$
in $\D \backslash \{ 0 \}$.
Further, for the remainder function
 $w(z)=(\log (1/|z|))^{-\beta}$ we see that
$$ w_z(z)=\frac{ \beta}{2} \, \frac{ 1}{z} \,  \frac{1}{(\log( 1/|z|))^{\beta+1}}\, .$$
\end{example}

\medskip

On the other hand, if $\kappa$ is H\"older continuous in $\D$, then
the potential--theoretic estimate (\ref{eq:integral}) combined with Lemma
\ref{lem:wachstum} immediately proves Theorem \ref{thm:exceptional} for the
case $\alpha=1$. So to finish the proof of Theorem \ref{thm:exceptional}
we finally need to prove Lemma \ref{lem:wachstum}.

\bigskip

{\bf Proof of Lemma \ref{lem:wachstum}.}

\nopagebreak
Without loss of generality we may suppose $\kappa(0)=-1$. Further, it is
advantageous to  take the regular conformal Riemannian metric
$\lambda(z)\, |dz| =e^{u(z)} \, |dz|=(|z| \log(1/|z|))^{-1} \sigma(z)\, |dz|$, where $0<
\liminf_{z \to 0} \sigma(z) \le \limsup_{z \to 0} \sigma(z)< \infty$ with curvature $\kappa$ into account.
In a first step we establish
\begin{equation}\label{eq:comp_kappa}
\displaystyle  C_1  \left( \log \frac{1}{|z|}
\right)^{-1} \le \lambda(z)\, |z| \log \frac{1}{|z|} -1 \le C_2 \left( \log \frac{1}{|z|}
\right)^{-1} 
\end{equation}
for all $z$ in some disk $K_{\varrho}$ and constants $C_1<0< C_2$ by
comparing  $\lambda(z)\, |dz|$ with suitable conformal Riemannian metrics. 

\medskip

For the moment we define for arbitrary $R>0$ on $K_R \backslash \{ 0 \}$ the
conformal Riemannian metrics
\begin{align*}
& \lambda_{R}^a(z)\, |dz|= \frac{\exp\left( a\,  (\log(R/|z|))^{-1}  \right)}{|z|
  \log(R/|z|)}\, |dz|\, ,\quad a>0\,, \\[-4mm]
\intertext{and} \\[-8mm]
&\lambda_R(z) \, |dz|= \frac{\exp\left((  \log(1/|z|) )^{-1}\right)}{|z|
  \log(R/|z|)} \, |dz|\, 
\end{align*}
with curvature
\begin{align*}
 &\kappa_{\lambda_R^a}(z) =-1 +2a^2 (\log(1/|z|))^{-2}+ O\left((\log(1/|z|))^{-3}\right)\\[-4mm]
\intertext{and}\\[-8mm]
&\kappa_{\lambda_R}(z)=-1 +2(1-2 \log{R})\,  (\log(1/|z|))^{-2} + 
O\left((\log(1/|z|))^{-3}\right) 
\end{align*}
respectively. We observe that
\begin{itemize}
\item[(i)] $\lambda_R^a(z)\to + \infty$ as $R \searrow |z|$ for fixed $|z|$;

\item[(ii)] $\lambda_R(z)\to 0$ as $R \to \infty$  for fixed $z$;
\item[(iii)] $R \mapsto \kappa_{\lambda_R^a}(z)$ and $R \mapsto \kappa_{\lambda_R}(z)$  are monotonically
decreasing. 
\end{itemize}
Now, choosing $a$ and $R$ appropriately will lead to (\ref{eq:comp_kappa}).

\medskip

For that we first find $C>0$ and $\varrho>0$ such that 
$$-C\le r(z)\le C \, ,\quad  z \in K_{\varrho} \, . $$

To derive the right inequality of (\ref{eq:comp_kappa}),  fix $a>0$ with $a^2
\ge C$. Then, shrinking $\varrho$ if necessary, we deduce that 
\begin{equation*}
\begin{split}
-\kappa(z)&= 1-r(z) \left( \log \frac{1}{|z|} \right)^{-2}
 \ge 
- \kappa_{\lambda_1^a}(z) +o((\log(1/|z|)^{-2}) + a^2 \left( \log
 \frac{1}{|z|} \right)^{-2}  \ge- \kappa_{\lambda_1^a}(z) 
\end{split}
\end{equation*}
for $z \in K_{\varrho}$.
 The monotonicity of $\kappa_{\lambda_R^a}$ implies now that for all $R$ with $\varrho
 < R <1$ 
$$-\kappa(z)\ge  - \kappa_{\lambda_R^a}(z)\, , \quad z \in
K_{\varrho}\, .$$
So $z \mapsto \log \lambda_R^a(z)$ is   a supersolution to
$\Delta u=-\kappa(z)\,  e^{2u}$  in $K_{\varrho} \backslash \{
0 \}$ for every $R \in (\varrho,1)$. Further,
by (i) we choose $R \in (\varrho,1)$ such that 
$\lambda(z) \le \lambda_R^a(z)$ for $|z|= \varrho$.
Therefore by Theorem \ref{thm:ext_max} we get
\begin{equation}\label{eq:lambda1}
 \lambda(z) \le \lambda_R^a(z) \quad \text{for}\, \,  z\in K_{\varrho}
\backslash \{ 0 \}\, 
\end{equation}
and finally
$$\lambda(z) \, |z| \, \log(1/|z|) \le \frac{\log(1/|z|)}{\log(R/|z|)}
\exp\left( a \, (\log(R/|z|))^{-1} \right)= 1+ (a -\log{R}) \, \frac{1}{\log(1/|z|)} +
\cdots\, .$$

\smallskip

The proof for the left inequality of (\ref{eq:comp_kappa}) runs similarly. Here we  choose
$R'>1$ such that $-(1-2 \log R') \ge C$. Then, again shrinking $\varrho$ if
necessary,  we have
\begin{equation*}
- \kappa(z)
 \le - \kappa_{\lambda_{R'}}(z) +o((\log(1/|z|))^{-2}) +(1-2 \log R') \left( \log
  \frac{1}{|z|}\right)^{-2}
\le - \kappa_{\lambda_{R'}}(z) 
\end{equation*}
 for $z \in K_{\varrho}$, which shows that for all $R>R'$
$$-\kappa(z) \le - \kappa_{\lambda_R}(z)\, ,  \quad  z
\in K_{\varrho}\, .$$
Thus $\log \lambda_{R}(z)$ is a subsolution to $\Delta u= -
\kappa(z) \, e^{2u}$  for every $R>R'$ and since $\lambda_R(z) \to 0$ for $R\to
\infty$ we can find an $R>R'$ such that
$\lambda(z) \ge \lambda_{R}(z)$ for $|z|=\varrho$.
Using again Theorem \ref{thm:ext_max} gives
\begin{equation}\label{eq:lambda2}
 \lambda(z) \ge \lambda_{R}(z) \quad \text{for} \, \,  z \in K_{\varrho}
 \backslash \{ 0 \}
\end{equation}
and consequently
$$ \lambda(z)\, |z| \, \log \frac{1}{|z|} \ge
\frac{\log(1/|z|)}{\log(R/|z|)}\, 
\exp\left( (\log (1/|z|))^{-1} \right)  =1+(1- \log R) \,  \frac{1}{\log
(1/|z|)} + \cdots\, . $$
This yields the left--hand side of inequality (\ref{eq:comp_kappa})
for $z \in K_{\varrho}$. 

\medskip

To see that (\ref{eq:q}) holds,
just note that $e^{w(z)}=\lambda(z)|z| \log(1/|z|)$ and
$$-\kappa(z) \; \! e^{2w(z)}-1=(e^{w(z)}-1) (e^{w(z)}+1) -\frac{r(z)
  e^{2w(z)}}{(\log(1/|z|))^2}\, .$$
Since $w(z)$ is bounded on every disk $K_R$, $R<1$, the result follows from
  (\ref{eq:comp_kappa}), shrinking $\varrho$ again, if necessary.\hfill{$\blacksquare$}

\subsection{Second derivative} \label{par:abl2}

We are finally left to establish the statements about the second derivatives  
of Theorem \ref{thm:hoelder}.
The proof is similar to the case of the
first derivatives above. Therefore, we restrict ourselves to indicate the 
necessary changes in the argument for these derivatives.

\bigskip

{\bf Proof of Theorem \ref{thm:hoelder}: Second derivatives.}\\
By assumption $\kappa$ is locally H\"older continuous in $\D$ and we may
assume a fixed H\"older exponent
$\gamma$, $0 < \gamma \le 1$, say.
Without loss of generality we may also assume $-a < \inf_{z \in \D} \kappa(z)<
\sup_{z \in \D} \kappa(z) <-A$, where $a$ and $A$ are positive constants and
$\kappa(0)=-1$.

\medskip

If $\alpha \le 0$ the result follows directly from the representation formula
(\ref{eq:representation1}) and Proposition
\ref{thm:potential}.

\medskip

For $0<\alpha <1$ 
we define $q(\xi)=-\kappa( \xi ) e^{2 v(\xi)}$ and note that $q$ is locally
H\"older continuous in $\D$ with exponent $\tilde{\gamma} =\gamma $ if
$0<\alpha <1/2$ and $\tilde{\gamma}=\min \{ \gamma, 2-2\alpha\}$ if $ 1/2 \le \alpha <1$.
Next, fix $R<1$, choose $z
\in K_{R/2}\backslash \{ 0 \}$ and set $r=|z|/2$.
Rearranging (\ref{eq:secondderivative}) yields

\begin{equation*}
\begin{split}
&\frac{\partial^2}{\partial x_l \, \partial x_j}
v(z)=\, \frac{\partial^2}{\partial x_l \, \partial x_j}h(z)  
-\frac{1}{2\pi} \,
\frac{q(z)}{|z|^{2\alpha}} \int \limits_{\partial K_r(z)}
  \frac{\partial}{\partial x_j} \log|z - \xi| \, n_{l} (\xi) \, |d\xi| \, +\\[2mm]
&\quad \frac{1}{2\pi}\iint \limits_{K_R\backslash K_r(z)} \frac{\partial^2}{\partial x_l \,
  \partial x_j}\log|z - \xi| \, \, \, \frac{q(\xi)}{|\xi|^{2\alpha}}\,
d\sigma_{\xi} 
+\frac{1}{2 \pi} \iint \limits_{ K_r(z)} \frac{\partial^2}{\partial x_l \,
  \partial x_j}\log|z - \xi| \,   \left( \frac{q(\xi)}{|\xi|^{2\alpha}}-
  \frac{q(z)}{|z|^{2\alpha} } \right)\, d\sigma_{\xi}  
\end{split}
\end{equation*}

for $l,j \in \{1,2\}$, where $h$ is harmonic on $K_R$.
Our aim is to show that the second derivatives of $v$ belong to  $O(|z|^{-2
  \alpha})$. For that we need only observe that
\begin{equation*}
\begin{split}
&\left|\frac{1}{2\pi}\iint \limits_{K_R\backslash K_r(z)} \frac{\partial^2}{\partial x_l \,
  \partial x_j}\log|z - \xi| \, \, \, \frac{q(\xi)}{|\xi|^{2\alpha}}\,
d\sigma_{\xi}\right| 
\le \sup_{\xi \in K_R} |q(\xi)|\, \, 
\frac{1}{2\pi}\iint \limits_{K_R\backslash K_r(z)} \frac{2}{|z - \xi|^2}
  \frac{1}{|\xi|^{2\alpha}}\,
d\sigma_{\xi}\le\frac{C_1}{|z|^{2\alpha}}
\end{split}
\end{equation*}
for some constant $C_1>0$  (see \cite[p.~215]{Els91}) and 
\begin{equation*}
\begin{split}
&\left| \frac{1}{2 \pi} \iint \limits_{ K_r(z)} \frac{\partial^2}{\partial x_l \,
  \partial x_j}\log|z - \xi| \, \,  \left( \frac{q(\xi)}{|\xi|^{2\alpha}}-
  \frac{q(z)}{|z|^{2\alpha} } \right)\, d\sigma_{\xi} \right| \le\\[1.5mm]
&\frac{1}{2\pi} \iint \limits_{K_r(z)} \frac{2}{|z - \xi|^2} \frac{|q(\xi)
  -q(z)|}{|\xi|^{2\alpha}}\, d\sigma_{\xi} + \frac{2M}{2\pi} \iint \limits_{K_r(z)}
  \frac{1}{|z- \xi|^2} \frac{(|\xi|^{\alpha} + |z|^{\alpha})\, \big|
  |\xi|^{\alpha} - |z|^{\alpha}\big|}{|z|^{2\alpha} |\xi|^{2\alpha}}\,
  d\sigma_{\xi}\le \frac{C_2}{|z|^{2\alpha}} \, , \\
\end{split}
\end{equation*}
where  $C_2$ is  some positive constant.

\smallskip

The case $\alpha=1$ runs similarly. At first pick $R<1/e^2$, define $q(\xi)=-
\kappa(\xi)\, e^{2w(\xi)}-1$ and put
  $M= \sup_{\xi \in K_R} |q(\xi)|$. Since $\kappa$ fulfills the hypotheses of Lemma \ref{lem:wachstum}
condition (\ref{eq:q}) holds for some disk $K_{\varrho}$.
Now, let $\tilde{\varrho}=\min \{R/2, \varrho \}$. Choose $z \in
K_{\tilde{\varrho}}$ and set $r=|z|/2$. Rewriting (\ref{eq:second_derivative-w}) gives
\begin{equation*}
\begin{split}
\frac{\partial^2}{\partial x_l \, \partial x_j}
w(z)=\, &\frac{\partial^2}{\partial x_l \, \partial x_j}h(z) + 
\frac{1}{2\pi}\iint \limits_{K_{\tilde{\varrho}}\backslash K_r(z)} \frac{\partial^2}{\partial x_l \,
  \partial x_j}\log|z - \xi| \, \, \, \frac{q(\xi)}{|\xi|^{2}(\log(1/|\xi|))^2}\,
d\sigma_{\xi} + \\[1.5mm]
&\frac{1}{2 \pi} \iint \limits_{ K_r(z)} \frac{\partial^2}{\partial x_l \,
  \partial x_j}\log|z - \xi| \, \,  \left( \frac{q(\xi)}{|\xi|^{2}(\log(1/|\xi|))^2}-
  \frac{q(z)}{|z|^{2}( \log(1/|z|))^2} \right)\, d\sigma_{\xi} - \\[1.5mm]
&\frac{1}{2\pi}\, \frac{q(z)}{|z|^{2}(\log(1/|z|))^2} \,  \int \limits_{\partial K_r(z)}
  \frac{\partial}{\partial x_j} \log|z - \xi| \, n_{l} (\xi) \, |d\xi|
\end{split}
\end{equation*}
for $l,j \in \{1,2\}$ and a harmonic function $h$ on $K_{\tilde{\varrho}}$.
To derive that the second derivatives of $w$ are in $O(|z|^{-2}
(\log(1/|z|))^{-2})$ as $|z| \to 0$,
we first note that
\begin{equation*}
\left|\frac{1}{2\pi}\iint \limits_{K_{\tilde{\varrho}}\backslash K_r(z)} \frac{\partial^2}{\partial x_l \,
  \partial x_j}\log|z - \xi| \, \, \, \frac{q(\xi)}{|\xi|^{2}(\log(1/|\xi|))^2}\,
d\sigma_{\xi}\right| \le
 \frac{C_3}{|z|^2 (\log(1/|z|))^2}\, , 
\end{equation*}
for some constant  $C_3$ where we have applied Lemma \ref{lem:wachstum}.
On the other hand by using the mean value theorem
and the H\"older continuity of $\kappa$ 
we find 
\begin{equation*}
\begin{split}
&\left| \frac{1}{2 \pi} \iint \limits_{ K_r(z)} \frac{\partial^2}{\partial x_l \,
  \partial x_j}\log|z - \xi| \, \,  \left( \frac{q(\xi)}{|\xi|^{2}(\log(1/|\xi|))^2}-
  \frac{q(z)}{|z|^{2}(\log(1/|z|))^2 } \right)\, d\sigma_{\xi} \right|
  \le\\[1.5mm]
&\frac{1}{2\pi} \iint \limits_{K_r(z)} \frac{2}{|z - \xi|^2} \frac{|\kappa(\xi)
  -\kappa(z)| e^{2w(\xi)}}{|\xi|^{2}(\log(1/|\xi|))^2}\, d\sigma_{\xi}+
\frac{1}{2\pi} \iint \limits_{K_r(z)} \frac{2}{|z - \xi|^2} \frac{|\kappa(z)|\,
  |e^{2w(\xi)} - e^{2w(z)}|}{|\xi|^{2}(\log(1/|\xi|))^2}\, d\sigma_{\xi}+\\[1.5mm]& \frac{1}{2\pi} \iint \limits_{K_r(z)}
  \frac{2}{|z- \xi|^2} |q(z)| \left| \frac{1}{|z|^2(\log(1/|z|))^2}  -
  \frac{1}{|\xi|^2(\log(1/|\xi|))^2}  \right| \, d \sigma_{\xi} \le
\frac{C_4}{|z|^2 (\log(1/|z|))^2} 
\end{split}
\end{equation*}
for some $C_4 >0$. Thus $w_{zz}(z)= w_{z\overline{z}}(z)=w_{\overline{z}\,
  \overline{z}}(z)=O\left(|z|^{-2}(\log(1/|z|))^{-2} \right)$ for $z\to 0$,  as desired. \hfill{$\blacksquare$}

\subsection{Proof of Theorem \ref{thm:2} and Theorem \ref{thm:3}} \label{par:thm23}

The proof of Theorem \ref{thm:2} is a straightforward application
of Theorem \ref{thm:hoelder} by noting that the metric $\lambda(z) \, |dz|$
has a representation of the form 
$$\lambda(z)=e^{u(z)}=\begin{cases} |z|^{-\alpha} e^{v(z)} & \text{ if } \alpha <1 \\[2mm]
   \displaystyle \frac{e^{w(z)}}{|z| \log (1/|z|)} & \text{ if } \alpha=1\,
   ,\end{cases}$$
where $u(z)$, $v(z)$ and $w(z)$ have the properties stated in Theorem
\ref{thm:hoelder} and Theorem \ref{thm:alpha=1,continuity}.
We note that the statements (a) and (b) of Theorem \ref{thm:2} remain
valid under the weaker assumption that the curvature $\kappa$ is only
continuous on $\D$. This follows from an inspection of the proofs in \S
\ref{par:conti}--\S \ref{par:abl2}. 

\medskip

In order to prove Theorem \ref{thm:3} we first note that
Theorem \ref{thm:hoelder} implies that
 a conformal Riemannian metric $\lambda(z) \, |dz|$ on
$\Omega$ with (H\"older continuous) curvature $\kappa : \Omega \cup \{ 0 \}
\to \R$ with $\kappa(0)<0$ is locally complete at $z=0$ if and only if 
$\log \lambda$ has order $\alpha=1$. In particular, the hyperbolic metric
$\lambda_{\Omega}(z) \, |dz|$ has also order $\alpha=1$.
Thus Theorem \ref{thm:3} follows immediately from Theorem \ref{thm:2}.

\section{Appendix: Potential Theory} \label{par:potential}

In this appendix we provide some essentially well--known but non--standard 
facts from potential theory such as a Poisson--Jensen formula and some
 differentiability properties of Newton's
Potential that are extensively and repeatedly used in the course of this
paper.

\begin{proposition}[Poisson--Jensen]\label{thm:poisson-jensen}
Let $u$ be a subharmonic function on $K_r$ such that
$u \in C^2(K_r \backslash \{ 0 \})$, $\Delta u \in L^1(K_r)$  and
$$\lim_{r\to 0} \frac{\sup_{|z|=r} u(z)}{\log(1/r)}=0\, .$$
Then 
$$u(z)=h(z) +\frac{1}{2\pi} \iint\limits_{K_r} \log|z - \xi| \, \Delta u\,
d\sigma_{\xi}\, , \qquad  z \in K_r\, ,$$ 
where $h$ is a harmonic function on $K_r$.
\end{proposition}

This can be deduced from Theorem 4.5.1 and Exercise 3.7.3 in \cite{Ransford}.

\begin{proposition}[Newton--Potential]\label{thm:potential}
\begin{itemize}
\item[(a)]
Let $r \le 1$ and $q:K_r \to \R$ be a bounded and integrable function. Then
for every $\alpha< 1$ the function
$$ 
\omega:K_r \to \R\, , \quad  z=x_1\!+i\!\,x_2\,  \mapsto \,  \frac{1}{2\pi} \iint
\limits_{K_r}\log|z- \xi| \, \, 
\frac{q(\xi)}{|\xi|^{2\alpha}} \, d\sigma_{\xi} \,  
$$
is continuous in $K_r$.\\
Further, $\omega \in  C^1(K_r)$  for $\alpha <1/2$ and $\omega \in C^1(K_r \backslash
\{ 0 \})$ for $1/2 \le \alpha <1$, where
$$
\frac{\partial }{\partial x_j}\omega(z)= \frac{1}{2\pi} \iint \limits_{K_r}
\frac{\partial}{\partial x_j} \log|z - \xi| \, \, \frac{q(\xi)}{|\xi|^{2\alpha}}
\, d\sigma_{\xi}
$$
for $j \in \{1,2\}$.\

\smallskip 
If, in addition, $q$ is locally H\"older continuous, then $\omega \in C^2(K_r)$ if $\alpha \le 0$  and   $\omega \in
C^2(K_r \backslash\{ 0 \})$ if $0 <\alpha<1$,
where
\begin{equation}\label{eq:secondderivative}
\begin{split}
\frac{\partial^2}{\partial x_l \partial x_j} \omega(z)=\frac{1}{2\pi }
\iint \limits_{K_3} \frac{\partial^2}{\partial x_l \partial x_j} \log|z-\xi|
\, \, \left(\frac{q(\xi)}{|\xi|^{2\alpha}}-\frac{q(z)}{|z|^{2\alpha}}
\right)\, d\sigma_{\xi}\\[2mm]
 -\,\frac{1}{2\pi} \frac{q(z)}{|z|^{2\alpha}} \int\limits_{\partial K_3}
\frac{\partial}{\partial x_j} \log|z- \xi|\,  n_l(\xi) \, |d\,\!\xi|\, .
\end{split}
\end{equation}
Here $(n_1(\xi),n_2(\xi))^T$ is the unit outward normal at the point $\xi \in
\partial K_3$ and $q$ is extended to vanish outside of $K_r$.
\item[(b)]
Let $r<1$ and let $q:K_r\to \R$ be a bounded and integrable function in
$K_r$.
Then the function
$$ \omega:K_r \to \R\, , \quad  z=x_1\!+i\!\,x_2\,  \mapsto \,  \frac{1}{2\pi} \iint
\limits_{K_r}\log|z- \xi| \, \, 
\frac{q(\xi)}{|\xi|^{2} (\log(1/|\xi|))^2} \, d\sigma_{\xi} \,$$
belongs to $C^1 (K_r \backslash \{ 0 \})$ and
$$
\frac{\partial }{\partial x_j} \omega(z)= \frac{1}{2\pi} \iint \limits_{K_r}
\frac{\partial}{\partial x_j} \log|z - \xi| \, \, \frac{q(\xi)}{|\xi|^{2}(\log(1/|\xi|))^2}
\, d\sigma_{\xi}
$$
for $j \in \{1,2\}$.

\smallskip

If, in addition, $q$ is locally H\"older continuous, then  $\omega \in C^2(K_r \backslash\{ 0 \})$. Further,
\begin{equation}\label{eq:second_derivative-w}
\begin{split}
\frac{\partial^2}{\partial x_l \partial x_j} \omega(z)=\frac{1}{2\pi }
\iint \limits_{K_3} \frac{\partial^2}{\partial x_l \partial x_j} \log|z-\xi|
\, \, \left(\frac{q(\xi)}{|\xi|^{2}(\log(1/|\xi|))^2}-\frac{q(z)}{|z|^{2}(\log(1/|z|))^2}
\right)\, d\sigma_{\xi}\\[2mm]
 -\,\frac{1}{2\pi} \frac{q(z)}{|z|^{2}(\log(1/|z|))^2} \int\limits_{\partial K_3}
\frac{\partial}{\partial x_j} \log|z- \xi|\,  n_l(\xi) \, |d\,\!\xi|\, , \hspace{1cm}
\end{split}
\end{equation}
where $(n_1(\xi),n_2(\xi))^T$ is the unit outward normal at the point $\xi \in
\partial K_3$ and $q$ is extended to vanish outside of $K_r$.
\end{itemize}
\end{proposition}

The result is standard for $\alpha \le 0$ (see \cite[p.54/55]{GT97}) and the
proof for $\alpha \le 0$ can be extended to the cases $0 <\alpha < 1$.


\end{document}